\documentclass{amsart}

\usepackage{graphicx}
\usepackage[nodisplayskipstretch]{setspace}
\usepackage[inline]{enumitem}

\usepackage{wrapfig}
\usepackage{subfig}
\usepackage{tikz}
\usepackage{amsfonts,amssymb,amsmath,amsthm}
\usepackage{float}
\usepackage{enumitem}
\usepackage{comment}
\usepackage{graphicx}
\usepackage[capitalise]{cleveref}
\usepackage{soul}
\usepackage{lineno}

\newtheorem{theorem}{Theorem}
\newtheorem{lemma}{Lemma}
\newtheorem{definition}{Definition}
\newtheorem{example}{Example}
\newtheorem{remark}{Remark}

\newtheorem{proposition}{Proposition}

\newcommand{\N}{\ensuremath{\mathbb{N}}}
\newcommand{\Naturals}{\ensuremath{\mathbb{N}}}
\newcommand{\Z}{\ensuremath{\mathbb{F}}}

\newcommand{\Q}{\ensuremath{\mathbb{Q}}}

\newcommand{\R}{\ensuremath{\mathbb{R}}}

\newcommand{\vset}{\mathcal{V}}
\newcommand{\height}{\texttt{height}}
\newcommand{\divides}[2]{{#1}|{#2} }
\newcommand{\notdivides}[2]{#1\not\hspace{2.5pt}\mid #2}
\newcommand{\idealRes}{\texttt{IdealRes}}
\newcommand{\monoRes}{\texttt{MonoRes}}

\newcommand{\degree} {r}

\newcommand{\poly}{\mathsf{p}}
\newcommand{\monomial}{\alpha}
\newcommand{\monomials}{\varset^{\oplus}}
\newcommand{\set}[1]{\ensuremath{\{ #1 \}}}
\newcommand{\suchthat}{\ensuremath{\, \vert \,}}
\newcommand{\aset}{A}

\newcommand{\cset}{\mathcal{C}}
\newcommand{\dset}{\mathcal{D}}

\newcommand{\sset}{S}
\newcommand{\fset}{\mathcal{F}}
\newcommand{\pset}{\mathcal{P}}
\newcommand{\Pset}{\mathbb{P}}
\newcommand{\gset}{\mathcal{G}}
\newcommand{\iset} {\mathcal{I}}

\newcommand{\barb}{\bar{b}}

\newcommand{\upclosure}[1]{#1{\uparrow}}

\newcommand{\monostrict}{\prec}
\newcommand{\monoorder}{\preccurlyeq}
\newcommand{\monolex}{\monoorder_{lex}}
\newcommand{\monodeg}{\monoorder_{deg}}

\newcommand{\termsof}{\mathcal{T}}
\newcommand{\htermsof}{\mathcal{HT}}
\newcommand{\hightermsof}[1]{\mathcal{HT}_{\!\!{#1}}}

\newcommand{\bigO}{\mathcal{O}}

\newcommand{\Ideal}[1]{\langle #1 \rangle}
\newcommand{\polynomials}{\Z_2[X]}
\newcommand{\ffun}{f}
\newcommand{\gfun}{g}

\newcommand{\nonresidual}{\mathcal{R}}
\newcommand{\mnonresidual}{\mathcal{R}_{\preccurlyeq}}
\newcommand{\Poly}{\mathsf{P}}
\newcommand{\F}{\mathbb{F}}

\newcommand{\mordering}{\triangleleft}
\newcommand{\pordering}{ \triangleleft}

\newcommand{\GB}{ \texttt{\bf Gb}}
\newcommand{\minel}{ {\tt Min}}
\newcommand{\EXP}{ e}
\newcommand{\degreeof}[2] {\texttt{Deg}_{#1}(#2)}

\newcommand{\wtdeg} {\texttt{wtDeg}}
\newcommand\wtmap {\texttt{u}}

\newcommand{\varset}{\mathbb{X}}

\newcommand{\barf} {\bar{f}}
\newcommand{\bars} {\bar{s}}

\newcommand{\barq} {\bar{q}}
\newcommand{\barc} {\bar{c}}

\newcommand{\barS} {\bar{S}}

\newcommand{\barQ} {\bar{Q}}
\newcommand{\barC} {\bar{C}}
\newcommand{\barB} {\bar{B}}
\newcommand{\barP} {\bar{\mathcal{P}}}
\newcommand{\p}{{\tt p}}

\begin{document}
	\title{Reduced Gr\"obner Bases With Double Exponential Cardinality  }
	
\author{Archana S Morye}
\author{Sreenanda S B}
\author{Prakash Saivasan}

\date{\today}

\address{University of Hyderabad, India }
\email{asmsm@uohyd.ac.in}
\address{University of Hyderabad, India }
\email{18mmpp02@uohyd.ac.in}
\address{The Institute of Mathematical Sciences, HBNI and CNRS IRL ReLaX,India}	
\email{prakashs@imsc.res.in}

	
\begin{abstract}
	In this article, we investigate the cardinality of Gr\"obner bases under various monomial orderings. We identify a family of polynomials $\fset$ and a criterion such that the reduced Gr\"obner basis $\GB_\monoorder(\fset)$ is double exponential in cardinality with respect to any monomial ordering $\monoorder$, which satisfies this criterion. We also show that the said criterion is satisfied by orderings such as the {lexicographic}, {degree lexicographic} and {weighted} orderings.
	
\end{abstract}

\maketitle

\section{Introduction}

In abstract algebra, most of the results are only existential but not constructive. For example, the Hilbert finite basis theorem \cite[p.76]{CLO1997}  states that 
ideals of polynomial rings are finitely generated. 
This theorem by itself does not give us a finite generating set for an ideal, but only tells that such a set exists.
 Sometimes it is not enough just to have any generating set.  In order to do computations efficiently, we need to have a generating set of particular kind.
 Such generating sets allow us to find solutions to  problems in the theory of polynomial rings,  the ideal description problem, the problem of solving system of polynomials and so on.

Gr\"obner basis is one such basis that was introduced by Bruno Buchberger in his PhD thesis \cite{Buc2006}.
The canonical version of such a basis is called the \emph{reduced Gr\"obner basis}.

Since the advent of Buchberger's algorithm, several new algorithms with improved running times have been developed \cite{Fau1999, Fau2002}.
In order to establish an optimal algorithm, the precise complexity of finding the reduced Gr\"obner basis has been widely explored. The best and worst case degree of Gr\"obner bases for a given set of polynomials has also been studied extensively \cite{Dub1990,Laz1983}. The growth rate of the coefficients while computing the Gr\"obner basis has also been well studied. For instance, \cite{Arn2000,Arn2003,Win1988} study polynomials with rational coefficients and provide an efficient algorithm for the same. In our work, we focus our attention on just the degrees and the cardinality of the reduced Gr\"obner bases.

 The first result on the complexity of Gr\"obner bases was due to  Mayr and Meyer.  In  \cite{MM1982,MR2013}, it was shown that the \emph{ideal membership} problem is {\sc Expspace-complete}. This also implies that the complexity of computing reduced Gr\"obner bases  can be large. The result was later improved by \cite{Yap1991} and \cite{BS1988}. In \cite{Yap1991}, it was shown that one can construct a set of polynomials such that any Gr\"obner basis for it has a polynomial with at least double exponential degree. In particular, the result showed that the double exponential degree can be achieved using a lesser number of polynomials in comparison to \cite{MM1982}. In \cite[Theorem II, p.199]{Huy1986}, Huynh constructed a family of polynomials such that under the degree lexicographic ordering, any Gr\"obner basis for it has a double exponential cardinality and degree.
This immediately establishes the fact that any algorithm that computes a reduced Gr\"obner basis (for a given set of polynomials) under the degree lexicographic ordering requires at least double exponential time.

In \cite{Kal2001}, Kalorkoti provides us with a set of polynomials for which the cardinality of the reduced Gr\"obner basis under degree reverse lexicographic ordering is small, whereas, the same under lexicographic ordering can be exponential. 
While this may be true in general, our work shows that there is a family of polynomials for which the cardinality can be uniformly double exponential in size and degree, in both the orderings.

In this work, we generalize the result in \cite{Huy1986} of Huynh. We show that there is a family of polynomials such that their reduced Gr\"obner basis has double exponential cardinality and degree under various monomial orderings. In fact, based on these polynomials, we develop a sufficient criterion such that any monomial ordering that satisfies this criterion is guaranteed to have at least double exponential cardinality for its reduced Gr\"obner basis.

\paragraph{Organisation and Overview}
In \cref{sec :prelims},  we introduce the  notations used in the paper. In particular, we define a {partial ordering}  $\pordering$ among the polynomials that will be crucially used throughout the paper.
In \cref{sec:upclose}, for any monomial ordering $\monoorder$, we characterize the set of head terms of the reduced Gr\"obner basis of a set $F$ as the minimal elements with respect to $\pordering$ of another set $\mnonresidual(F)$. This allows us to reason about the cardinality of the reduced Gr\"obner basis in terms of the number of minimal elements of $\mnonresidual(F)$ with respect to $\pordering$.

In \cref{sec:SEGB}, we construct a family of polynomials and establish a criterion for a monomial ordering that can produce a double exponential  sized reduced Gr\"obner basis. 
 In \cref{sec:MonoOrder}, we show that the reduced Gr\"obner basis for this set of polynomials under the lexicographic ordering, degree lexicographic ordering and weighted ordering can be of double exponential size and degree.  For this, we simply show that each of these orderings satisfies the identified criterion.

We would  like to emphasise  that we prove the results by directly considering system of polynomials instead of  {semi-Thue system} that was used in \cite{Huy1986,MM1982}.  We believe this will provide a direct intuition on what parameters of the  system of polynomials contribute to  the double exponential  blow-up.

\section{Preliminaries}\label{sec :prelims}

We use $\Naturals$ to represent the set of all natural numbers. 
 For any $m<n \in \N$, we let $[m,n]=\set{m,\ldots,n}$.  
We use  $f(n) = \bigO(g(n))$ to mean $f(n)$ is asymptotically bounded above by $g(n)$, that is,  there exist $n_0, c \in \Naturals$ such that for all  $n \geq n_0$, we have $f(n) \leq c g(n)$.

Let $X=\set{x_1,\ldots,x_n}$, then $X^\oplus$ denotes the collection of all monomials (sometimes referred to as terms)  in $x_1,\ldots,x_n$.  
We denote $\Z_2[X]$ as the polynomial ring $\Z_2[x_1,\ldots,x_n]$, where  $\Z_2$ denotes the field  with two elements $\{0,1\}$.
For any monomial $t = x_1^{\degree_1}\cdots x_n^{\degree_n}$, we will use $\degreeof{}{t} = \degree_1 + \cdots +\degree_n$ to refer to the \emph{degree of the monomial}.
For any polynomial $f$, we will use $\termsof(f)$ to denote the set of all the terms (monomials) appearing in $f$, and $\degreeof{}{f} = \text{Max}(\{ \degreeof{}{t} \mid t \in \termsof(f) \})$.

\begin{definition}[Monomial Ordering $\monoorder$]
 A total ordering $\monoorder$ on $X^\oplus$ is called a \emph{monomial ordering} if it satisfies the following conditions:
\begin{enumerate}
\item $1\monoorder t$ for all $t\in X^\oplus$.
\item For any $t_1,t_2,t_3 \in X^\oplus$, if $t_1 \monoorder t_2 $, then $t_1  t_3 \monoorder t_2  t_3$. 
\end{enumerate}
We will use $t_1 \monostrict t_2$ to denote $t_1 \monoorder t_2$ and $t_1 \neq t_2$.
	 \end{definition}
 	 We let  $\hightermsof{\monoorder}$ be the  function that returns   the  highest term of the polynomial with respect to $\monoorder$ that is, it returns the term $t\in\termsof(f)$ such that $s\monoorder t$ for all $s\in \termsof(f)$. 
We extend the definition of $\monoorder$ to polynomials as follows.
	 \begin{definition}[Monomial Ordering $\monoorder$ for Polynomials]
	For any two polynomials $f$ and $g$, we say $g \monoorder f$ if and only if $g = f$ or one of the following holds: Either	$\hightermsof{\monoorder}(g) \monostrict \hightermsof{\monoorder}(f)$, or $\hightermsof{\monoorder}(g)=\hightermsof{\monoorder}(f)$ and $g-\hightermsof{\monoorder}(g) \monoorder f-\hightermsof{\monoorder}(f)$ inductively. 
	\end{definition}

For any  $F \subseteq \Z_2[X]$, we let $\hightermsof{\monoorder}(F)=\set{\hightermsof{\monoorder}(f) \suchthat f\in F}$. Whenever the ordering is clear from the context, we will omit the reference to ${\monoorder}$ in $\hightermsof{\monoorder}$ and simply use $\htermsof$.

\begin{example}
\label{ex:monoorder} 
Let $t=x^{\degree_1} y^{\degree_2}$ and $s=x^{\degree'_1} y^{\degree'_2}$ in $\Z_2[x,y]$. Then,
\begin{enumerate}
\item  The lexicographic ordering is given as follows: $t\monoorder_{lex} s$ if and only if $\degree_1<\degree'_1$ or $\degree_1=\degree'_1$ and $\degree_2\leq \degree'_2$.

\item The degree lexicographic ordering is given as follows: $t\monoorder_{deg} s$ if and only if $\degree_1+\degree_2 < \degree'_1+\degree'_2$ or if $\degree_1+\degree_2 =\degree'_1+\degree'_2$, $t\monoorder_{lex} s$.  

\item  Let $\wtmap:\set{x,y}\rightarrow \R$ be given by $\wtmap(x)=\sqrt{2}$ and $\wtmap(y)=\sqrt{3}$. Then $t\monoorder_{\wtmap}s$ if and only if $\degree_1\sqrt{2}+\degree_2\sqrt{3}\leq\degree'_1\sqrt{2}+\degree'_2\sqrt{3}$. 
Note that since $\sqrt{2}$ and $\sqrt{3}$ are linearly independent over $\Q$, if $\degree_1\sqrt{2}+\degree_2\sqrt{3}=\degree'_1\sqrt{2}+\degree'_2\sqrt{3}$ then $t=s$, which makes it a total order. This ordering is an example of the weighted ordering.
\end{enumerate}
\end{example}

 Recall that an \emph{ideal} is a subset of polynomials $\iset \subseteq \Z_2[X]$ such that $\iset$ is an additive group and
for any $f \in  \Z_2[X]$ and $g \in \iset$, $f g \in \iset$.  Let $F$ be a subset of $\polynomials$, then $\Ideal{ F} $ denotes the ideal generated by $F$, that is, $\Ideal{F}$ is the smallest ideal with respect to inclusion containing $F$ in $\Z_2[X]$.  

\begin{definition}
Given an ideal $\iset$, we say a finite set of polynomials $\gset \subset \iset$ is a \emph{Gr\"obner basis} with respect to a monomial ordering $\monoorder$, if $\Ideal{\hightermsof\monoorder(\iset)} = \Ideal{\hightermsof\monoorder(\gset)}$.
We say $\gset$ is a \emph{reduced Gr\"obner basis} if it is a Gr\"obner basis, and for all $f \in \gset$ and $t\in\termsof(f), t \notin \Ideal{\hightermsof\monoorder(\gset \setminus \{\ffun\})}$.  For any set of polynomials $F$, we let $\GB_{\monoorder}(F)$ denote the reduced Gr\"obner basis of $\Ideal{F}$. Whenever the ordering is clear from the context, we will use $\GB(F)$ instead.
\end{definition}

For any two polynomials $\ffun,\gfun \in \polynomials$, we say $\ffun \xrightarrow[F]{} \gfun$ if and only if $\ffun + \gfun \in \Ideal{F}$. We define $\idealRes_F(\ffun) = \{g \mid \ffun \xrightarrow[F]{} \gfun \}$. For any subset $S\subseteq \Z_2[X]$, we denote $\idealRes_F(S)=\bigcup_{s\in S} \idealRes_F(s)$. 
We use $\monoRes_{F}(\ffun)$ to refer to  $\idealRes_{F}(\ffun)\cap X^\oplus$.

\begin{example}
\label{ex:idealres}
Consider the set $F=\set{f_1,f_2}\subseteq \Z_2[x,y]$ where $f_1=xy+1$ and $f_2=x^2+y$. Let $f=x^3y+x$ and $g=x+y$, then $f+g=x^3y+y=x^2f_1+f_2$. Hence we have $g\in \idealRes_F(\ffun)$. 
\end{example}

We now introduce an ordering $\mordering$ over the monomials  $ X^\oplus$ and also extend it over $ \polynomials$. This will be crucial to characterize the ideal generated by the head terms of the reduced Gr\"obner basis of a given set of polynomials.

\begin{definition}[{  $\mordering$ ordering over $ X^\oplus$}]

For any two monomials  $t = x_1^{\degree_1}\cdots x_n^{\degree_n} \in X^\oplus$ and  $s = x_1^{\degree'_1}\cdots x_n^{\degree'_n} \in X^\oplus$, we say $t \mordering s$ if and only if for all $i \in [1,n]$, $\degree_i \leq \degree'_i$. 
It is easy to see that $\mordering$ induces a partial ordering on $X^\oplus$. We extend this ordering to polynomials as follows.
\end{definition}

\begin{definition}[{   $\pordering$ ordering over $ \polynomials$}]
Let $p$ and $q$ be any two polynomials in $\polynomials$. We say $p\pordering q$ if and only if   there is an injective mapping $\ffun: \termsof(p) \to \termsof(q) $ such that $t  \mordering \ffun(t)$ for every $t\in \termsof(p)$.  The relation $\pordering$ induces a partial ordering on $\polynomials$. 
\end{definition}

\begin{example}\label{ex:triangleorder}
{Let $p=x^3y+xy^2$ and $q=x^4+x^3yz+x^2y^2z$ be two polynomials in $\Z_2[x,y,z]$. Define a map $\ffun: \termsof(p) \to \termsof(q) $ with $x^3y\mapsto x^3yz$ and $xy^2\mapsto x^2y^2z$. Here, $x^3y \pordering x^3yz$ and $xy^2 \pordering x^2y^2z$. Hence $p\pordering q$.} 

Consider another polynomial $r=x^3y+xy^2+z\in \Z_2[x,y,z]$. Then $r$ and $q$ are incomparable  because we  cannot define an injective map from $\termsof(r)$ to $\termsof(q)$ as well as from  $\termsof(q)$ to $\termsof(r)$ that satisfies the above property. 
\end{example}

For any set of polynomials $F\subseteq \Z_2[X]$, we define the minimum of $F$ with respect to $\pordering$ as $\minel_{\pordering}(F)=\set{p\in F| \text{ for all } q\in F \text{ if } q\pordering p, \text{ then } q=p}$.

\begin{definition} Given  a subset $\aset \subseteq \polynomials$, we define its \emph{upward closure (denoted $\upclosure{\aset}$)} as 
$\upclosure{\aset} = \{ s  \in\ \polynomials \mid \text{ there is a  } p \in \aset\text{ such that } p \pordering s \}.$
We say a set $\aset$ is  \emph{upward closed} if  $\upclosure{\aset} = \aset$.
 \end{definition}

\begin{example}
\label{ex:upclosure}
Let $\aset= \set{x^2y,y^3}$ be a subset of $\Z_2[x,y]$. Then $\upclosure{\aset} = \set{x^2yp_1+p_2,y^3p_1+p_2\mid p_1,p_2\in \Z_2[x,y]\text{ and }p_1\neq 0}$. 
\end{example}

\section{Characterising the highest terms of a Gr\"obner basis}\label{sec:upclose}

Let $\monoorder$ be any monomial ordering, we will show that the highest terms of a Gr\"obner basis with respect to  $\monoorder$ can be characterized as the minimal elements of an upward closed set with respect to $\pordering$.
Let $F \subseteq \polynomials$ be any set of polynomials, we define its \emph{reducible set} of polynomials $\mnonresidual(F)$ (simply $\nonresidual(F)$ if $\monoorder$ is clear from the context) as  follows:
$$ \mnonresidual(F)=\set{f\in \Z_2[X]\mid \text{ there is a  } g\prec f \text{ such that} f+g\in \Ideal{F}}$$
We show that $\nonresidual(F)$  is precisely the upward closure of the set of highest terms of the reduced Gr\"obner basis $\GB(F)$. Since the size of a reduced Gr\"obner basis is at least as much as the number of highest terms, it is at least as much as the number of minimal elements in $ \nonresidual(F)$.

We first show in the following lemma 
that the highest term of any polynomial in $\Ideal{F}$ is contained in $\nonresidual(F)$ and also show that $\nonresidual(F)$ is upward closed with respect to $\pordering$.

\begin{lemma}\label{lem:upclose:nonres}
	Let $\nonresidual(F)$ be the set defined as above. If $\gfun\in \Ideal{F}$, then $\htermsof(\gfun) \in\nonresidual(F) $. Furthermore, $\nonresidual(F)$ is upward closed with respect to the order $\pordering$.
\end{lemma}
\begin{proof}
 Given that $g\in \Ideal{F}$, we can write $\gfun= \htermsof(\gfun) + \gfun' $. Clearly, $\gfun'\monostrict \htermsof(\gfun) $. This gives  $\htermsof(\gfun)\xrightarrow[F]{} \gfun'$, hence $\htermsof(\gfun) \in\nonresidual(F) $. This proves the first statement.

 For any  $\ffun \in \polynomials$  and $\gfun \pordering \ffun$ with $\gfun \in \nonresidual(F)$, to prove that $\nonresidual(F)$ is upward closed, we need to find a $p \monostrict\ffun$ such that $\ffun \xrightarrow[F]{}p $ (i.e, $\ffun \in \nonresidual(F)$).
Since $\gfun \in \nonresidual(F)$, by definition we have $\gfun \xrightarrow[F]{} q$ for some  $q \monostrict \gfun$.	Let $T=\termsof(\gfun)\cap \termsof(q)$. 
 
 \textbf{Case $T=\emptyset$:\quad}
 Write $\gfun=\htermsof(g)+\gfun'$, so $\gfun+q=\htermsof(\gfun)+\gfun'+q$. It is easy to see that $\gfun'+q \monostrict \htermsof( \gfun)$. 
By our assumption, we have $\gfun+q  \in \Ideal{F}$, hence $ \htermsof(\gfun)+ \gfun'+q  \in \Ideal{F}$. Since $\gfun \pordering \ffun$, there is a term $t \in \termsof(\ffun)$ such that   $\htermsof(\gfun) \pordering t$, hence $t = \alpha\htermsof(\gfun)$ for some $\alpha \in X^\oplus$. Now, we have $ \alpha (\htermsof(g)+\gfun'+q) = t+ \alpha (\gfun'+q)\in \Ideal{F} $.    
Let $\ffun =t + \ffun'$, which implies $\ffun + \ffun'+  \alpha (\gfun'+q)  \in \Ideal{F}$, then $\ffun \xrightarrow[F]{}  \ffun'+  \alpha (\gfun'+q)$.   
Now all we need to show is that $ \ffun'+  \alpha (\gfun'+q) \monostrict \ffun$. We already know that $\ffun' \monostrict \ffun $ so it is sufficient to show $\alpha (\gfun'+q) \monostrict f$. But notice that $ \alpha (\gfun'+q) \monostrict \alpha \htermsof(g) = t \monoorder f$.

 \textbf{Case $T\neq\emptyset$:\quad}
If $T\neq \emptyset$, we let $h=\sum_{\beta\in T} \beta$. Let $g''=h+g$ and $q''=h+q$, then $g+q =g''+q''  \in  \Ideal{F}$,  $q'' \monostrict g'' $, and $g'' \pordering g $. Hence $\gfun'' \in \nonresidual(F)$,  $g'' \pordering  f$, and $\termsof(\gfun'')\cap \termsof(q'') =\emptyset$. Now, replacing $g$ and $q$ by $g''$ and $q''$ in the above case completes the proof.
\end{proof}

Lemma \ref{lem:sf2gb} provides a connection between $\nonresidual(F)$ and $\GB(F)$.  In particular, it proves that the set $\nonresidual(F)$ is the upward closure of the set of highest terms in the reduced Gr\"obner basis of $\Ideal{F}$. 

\begin{lemma}\label{lem:sf2gb}
Let 	$\nonresidual(F)$ be the set defined as above. Then,  
		$ \nonresidual(F) = \upclosure{ \htermsof(\GB(F))}.$
\end{lemma}
\begin{proof}
To prove $ \upclosure{ \htermsof(\GB(F))} \subseteq \nonresidual(F)$, we first observe that $\GB(F) \subseteq \Ideal{F} $,
hence by the first statement of Lemma \ref{lem:upclose:nonres}, $\htermsof(\GB(F)) \subseteq \nonresidual(F)$. The result follows as we have already shown in Lemma \ref{lem:upclose:nonres} that $\nonresidual(F)$ is upward closed.

For the other direction, we show that $ \nonresidual(F) \subseteq \upclosure{ \htermsof(\GB(F))}$.
Let $\ffun \in \nonresidual(F)$, and for any term $t \in \termsof(\ffun)$, we have $t \pordering \ffun$. Therefore, it is enough to show that  for some $t \in \termsof(\ffun)$,
$ t \in \upclosure{ \htermsof(\GB(F))}$.

Since $\ffun \in \nonresidual(F)$, there is a $\ffun' \prec \ffun$ such that $\ffun + \ffun' \in \Ideal{F}$. This means $\ffun + \ffun' = \Sigma_{i=1}^n a_i \cdot g_i$, where $\GB(F) = \{g_1,\ldots,g_n\}$. Further, we have for all $i\in [1,n]$, $ \htermsof(a_i \cdot g_i) \monoorder \htermsof(\ffun + \ffun') $ \cite[p.32]{AL1994}. From this, we obtain an $i_0 \in [1,n]$ such that $ \htermsof(a_{i_0} \cdot g_{i_0}) = \htermsof(\ffun + \ffun') $.  Since $\ffun' \prec \ffun$, we get $\htermsof(\ffun + \ffun') \in  \termsof(\ffun)$ and then the required $t$ is $\htermsof(\ffun + \ffun')= \htermsof(a_{i_0})\cdot  \htermsof( g_{i_0})$. This implies $\htermsof( g_{i_0})\pordering t$ and hence  $t \in \upclosure{ \htermsof(\GB(F))}$.
\end{proof}

As a consequence of the above lemma, we obtain the following theorem which states that the number of minimal elements of $\nonresidual(F)$ is equal to the cardinality of the reduced Gr\"obner basis.
\begin{theorem}\label{thm:cardGB}
	For any set of polynomials $F$, $ |\GB(F)| = |\minel_\pordering(\nonresidual(F)) |.$
\end{theorem}
\begin{proof}
Since $\GB(F)$ is a	 reduced Gr\"obner basis, for any $\ffun_1, \ffun_2 \in \GB(F)$, $\htermsof(\ffun_1) $ and $\htermsof(\ffun_2)$ are not comparable with respect to $\pordering$. This implies that $\minel_\pordering( \htermsof(\GB(F))) = \htermsof(\GB(F))$. From this and \cref{lem:sf2gb},
 we get
\[ 
 \minel_\pordering( \nonresidual(F)) = \minel_\pordering( \upclosure{ \htermsof(\GB(F))})= \minel_\pordering({ \htermsof(\GB(F))})=\htermsof(\GB(F)).
\]
  Since $|\htermsof(\GB(F))|=|\GB(F)|$, we obtain $ |\GB(F)| = |\minel_\pordering(\nonresidual(F)) |$.
\end{proof}

\begin{example}
\label{ex:residual}
 Consider the set $F=\set{xy+1,x^2+y}$ with respect to lexicographic ordering with $y\preccurlyeq_{lex} x$. Its reduced Gr\"obner basis is, $\GB(F)=\set{x+y^2,y^3+1}$ and hence, $ \htermsof(\GB(F))=\set{x,y^3}$. Now, consider the set $A=\set{xp_1+p_2,y^3p_1+p_2\mid p_1,p_2\in \Z_2[x,y] \text{ where }p_1\neq 0}$. Observe that any polynomial in this set belongs to $\nonresidual(F)$, since:

\begin{itemize}
\item $xp_1+p_2\xrightarrow{x+y^2} y^2p_1+p_2$, and $y^2p_1+p_2\monolex xp_1+p_2$. 
\item $y^3p_1+p_2\xrightarrow{y^3+1}p_1+p_2$, and $p_1+p_2 \monolex y^3p_1+p_2$. 
\end{itemize}
Now, $\Z_2[X]\setminus A=\set{y^2,y^2+y,y^2+y+1,y^2+1,y,y+1,1}$, and it is easy to verify that none of these belong to $\nonresidual(F)$. Hence $\nonresidual(F)=A=\upclosure{ \htermsof(\GB(F))}$. Moreover, $\minel_\pordering( \nonresidual(F))=\set{x,y^3}=\htermsof(\GB(F))$.  
\end{example}

\section{A criterion to obtain double exponential Gr\"obner bases} \label{sec:SEGB}

In this section, for each $n \in \Naturals$ we will construct a set of polynomials $\fset(n)$ with $|\fset(n)|=\bigO(n)$ and provide a criterion such that any monomial ordering $\monoorder$ that satisfies this criterion has  $|\GB_\monoorder(\fset(n))| \geq 2^{2^n}$.  
In fact, we show that there are at least $2^{2^n}$ polynomials each with a double exponential degree  in $\GB_\monoorder(\fset(n))$.
For the sake of readability, we will let $\EXP(n) = 2^{2^n}$. Any notation defined in this section will be used throughout the rest of the paper.

In order to construct the set of polynomials $\fset(n)$, we crucially use the set of polynomials that was used in \cite[p.316]{MM1982}.
We will augment this set of polynomials with additional polynomials and show that the reduced Gr\"obner basis of it has $\EXP(n)$ cardinality. We will fix $n \in \Naturals$, and show how to obtain $ \fset(n)$. In fact, we will simply refer to $\fset(n)$ as $\fset$ from here onwards. The variables that we use in the set will be as follows.

\begin{equation*}
	\begin{split}
		\varset  =\  &  V \cup  \bar{V}    \cup \{ s,\ell,c,\barc,b,\barb\}, \text{ where }  
		V =  \bigcup_{i=0}^n V_i,  \quad   \bar{V} =  \bigcup_{i=0}^n \bar{V}_i  \text{ and }\\
		& V_i = \{s_i,f_i,q_{1i},q_{2i},q_{3i},q_{4i},c_{1i},c_{2i},c_{3i},  c_{4i},b_{1i},b_{2i},b_{3i},b_{4i}  \}	\\
		& \bar{V}_i = \{\bars_i,\barf_i,\barq_{1i},\barq_{2i},\barq_{3i},\barq_{4i},\barc_{1i},\barc_{2i},\barc_{3i},  \barc_{4i},\barb_{1i},\barb_{2i},\barb_{3i},\barb_{4i}  \}	.\\
	\end{split}
\end{equation*} 

We now provide the required  set of polynomials, $\fset$ as a union of three sets of polynomials that is, $\fset = \pset \cup \barP \cup \gset$. We will define the set $\pset$  of polynomials as $\bigcup_{m=0}^n\pset_m $, where each $\pset_m$ is defined as follows.

\begin{eqnarray}
	\pset_0 = \{ & & b_{i0}^2 c_{i0} f_0+ c_{i0} s_0 \mid i \in [1,4] \}  \label{pm:1}\\
	\pset_m = \big\{& & q_{1m}c_{1(m-1)}s_{m-1}+s_m  \label{pm:2}\\ 
	                &       &    q_{2m}c_{2(m-1)}s_{m-1}+q_{1m} b_{1(m-1)}c_{1(m-1)}f_{m-1}\label{pm:3}\\
	                &  &  q_{3m}c_{3(m-1)}f_{m-1}+ q_{2m} c_{2(m-1)} f_{m-1}  \label{pm:4}\\
	                &  & q_{3m} b_{1(m-1)}c_{3(m-1)}s_{m-1}+ q_{2m} b_{4(m-1)}c_{2(m-1)}s_{m-1}\label{pm:5}\\     
	                & & q_{4m}b_{4(m-1)}c_{4(m-1)}f_{m-1}+q_{3m} c_{3(m-1)} s_{m-1} \label{pm:6}\\ 
	                & &q_{4m}c_{4(m-1)}s_{m-1}+ f_m \}  \ \ \ \       \label{pm:7}  \cup \\ 
	            & & \{q_{2m}b_{3(m-1)}b_{im}c_{im}f_{m-1}+ q_{2m}b_{2(m-1)}c_{im}f_{m-1}\mid i \in [1,4]\}      \label{pm:8}
          \end{eqnarray}

We define $\barP$ very similar to $\pset$ by replacing all the normal variables with their barred counterpart. That is, it is defined as $\bigcup_{m=0}^n\barP_m$, where each  $\barP_m$ is defined as follows.

\begin{eqnarray}
   	\barP_0 = \{ & &  \barb_{i0}^2 \barc_{i0} \barf_0 + \barc_{i0} \bars_0 \mid i \in [1,4] \}\\
   	\barP_m = \big\{& &   \barq_{1m} \barc_{1(m-1)} \bars_{m-1}+\bars_m\\
   	&       &   \barq_{2m} \barc_{2(m-1)} \bars_{m-1}+\barq_{1m} \barb_{1(m-1)} \barc_{1(m-1)}\barf_{m-1}\\
   	&  & \barq_{3m}\barc_{3(m-1)}\barf_{m-1} +{ \barq_{2m}}\barc_{2(m-1)}\barf_{m-1}  \\
   	&  & \barq_{3m} \barb_{1(m-1)}\barc_{3(m-1)}\bars_{m-1}+ \barq_{2m}\barb_{4(m-1)}\barc_{2(m-1)}\bars_{m-1}\\    
   	& & \barq_{4m}\barb_{4(m-1)}\barc_{4(m-1)}\barf_{m-1}+\barq_{3m}\barc_{3(m-1)}\bars_{m-1}\\
   	& & \barq_{4m}\barc_{4(m-1)}\bars_{m-1}+ \barf_m \}  \ \ \ \        \cup \\
   	& &\{\barq_{2m}\barb_{3(m-1)}\barb_{im}\barc_{im} \barf_{m-1}+\barq_{2m}\barb_{2(m-1)} \barc_{im} \barf_{m-1}\mid i \in [1,4] \}    
   	 \}
   \end{eqnarray}

We define the set $\gset$ of polynomials as,
     \begin{equation}
        	\begin{split}
        	 \gset = \ 	\{ b_{4n}\ell b+\ell c,\ &  b_{4n}\ell \barb + \ell \barc,\ c_{4n} f_n+ \ell,\     \barc_{4n}\barf_n +c_{4n}s_n, \\
        	 			&  \barb_{4n}c_{4n}s_n +c_{4n}s_n b,\ \barb_{4n}c_{4n} s_n +c_{4n}s_n\barb,\  \barc_{4n} \bars_n+ s \}.
         	\end{split}
        \end{equation}
        
Note that $|\fset(n)| = 20n+15$, hence the size of $\fset(n)$ is $\bigO(n)$.  We now identify the criterion that will provide us with a scheme to obtain a Gr\"obner basis of high cardinality for any monomial ordering $\monoorder$ that satisfies it. 

In order to show that $|\GB(\fset)|\geq \EXP(n)$, by Theorem \ref{thm:cardGB} it is sufficient  to prove that $|\minel_{\pordering} (\nonresidual(\fset))|\geq \EXP(n)$. Our  strategy is to construct a set $\cset$ with cardinality at least $\EXP(n)$ and establish a criterion    that will ensure  $\cset\subseteq \minel_{\pordering} (\nonresidual(\fset)) $. The criterion will involves two conditions on any monomial ordering.
The  set $\cset$ is defined as follows
\begin{equation}
	\cset = \{\ell \barc^{m_1} c^{m_2}  \mid m_1 +m_2 = \EXP(n) \}
\end{equation}
Note that the cardinality of $\cset$ is  $\EXP(n)+1$.  In order to show $\cset\subseteq \minel_{\pordering} (\nonresidual(\fset)) $, we need $\cset \subseteq \nonresidual(F)$ and there is no polynomial  (in fact monomial since  $\cset$ contains only monomials) $\alpha' \in \nonresidual(\fset)$  such that $\alpha' \neq \alpha$ and $\alpha' \pordering \alpha$ for any $\alpha \in \cset$.  

\begin{figure}[H]%
	\centering
	\subfloat[ $\cset \subseteq \nonresidual(\fset)$]{{\includegraphics[width=5.5cm]{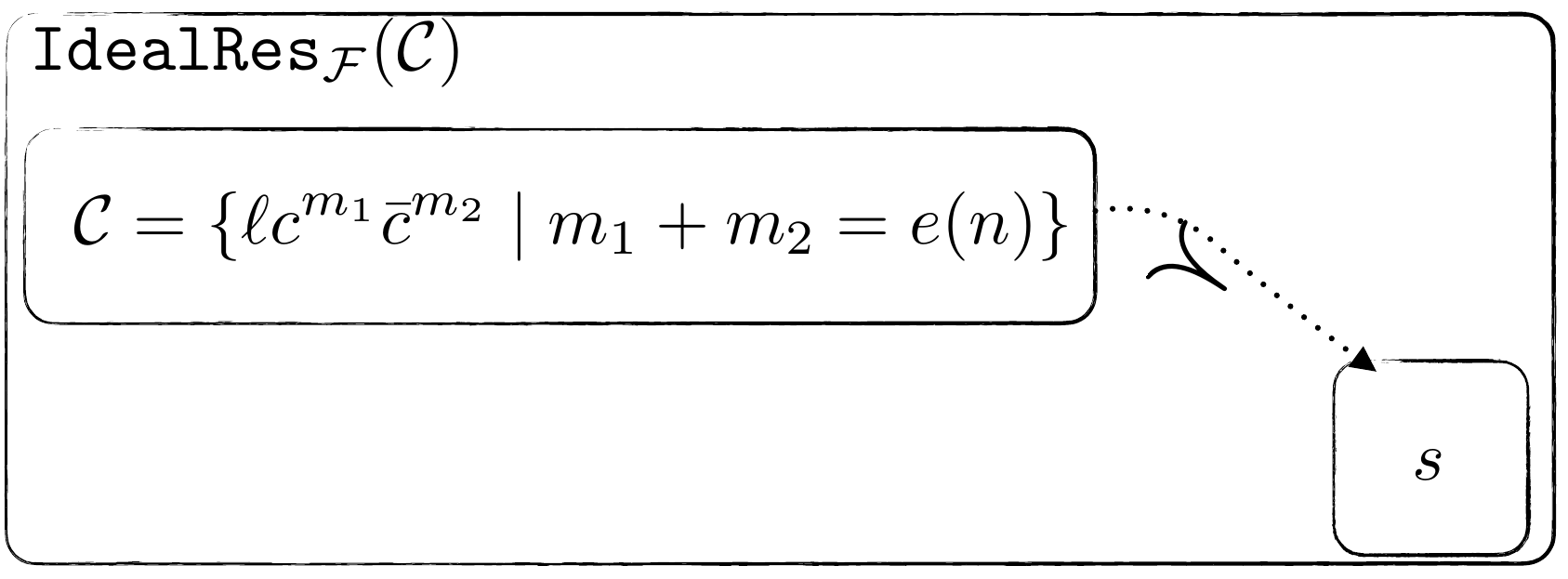} }} \label{fig:one}%
	\qquad
	\subfloat[ $\dset \cap \nonresidual(\fset)=\emptyset$]{{\includegraphics[width=5.5cm]{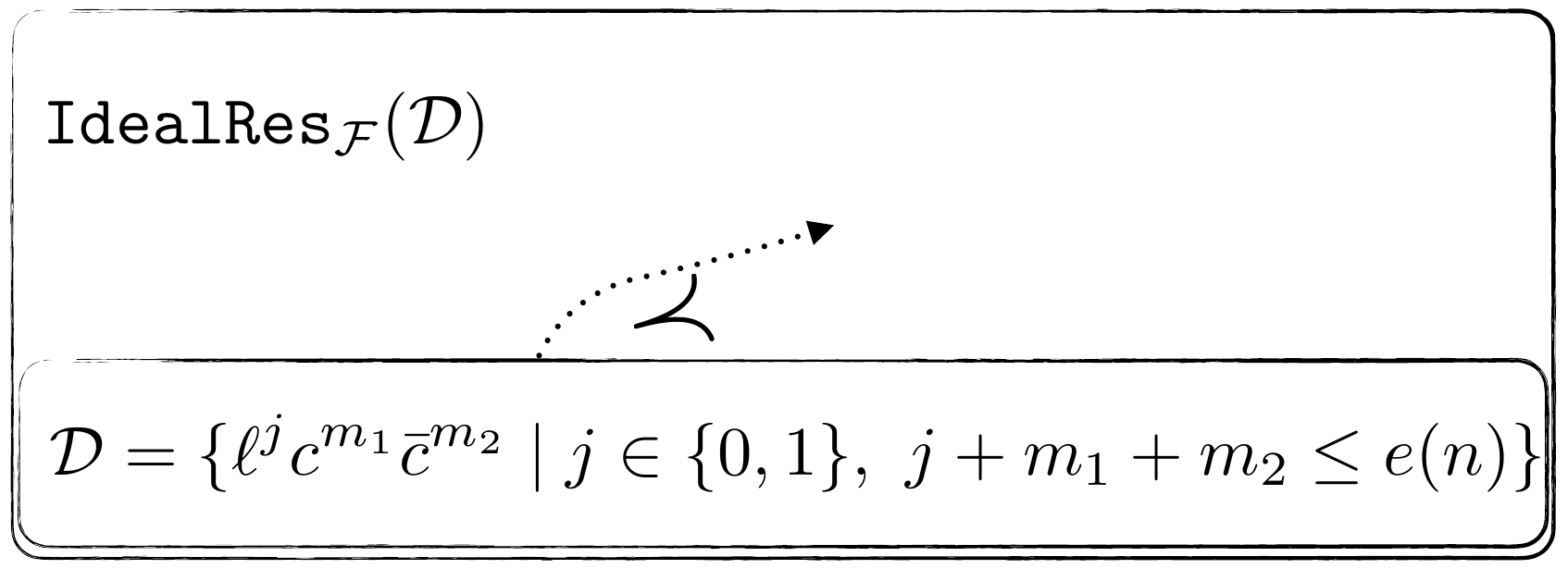} }}\label{fig:two}%
	\caption{Illustration for Theorem \ref{thm:main}}%
	\label{fig:example}%
\end{figure}

Corresponding to each $\alpha \in \cset$, if  there is a polynomial $f \in \idealRes_{\fset}(\alpha)$ such that $f \monostrict \alpha$, then by definition, $\cset \subseteq \nonresidual(\fset)$. 
The polynomials given in $\fset$ are such that $s \in \idealRes_{\fset}(\alpha) $ {(this, we will prove in Proposition \ref{prop:subseteq})}. If the monomial ordering is such that, $s \monostrict \alpha$ then we will get $\cset \subseteq \nonresidual(\fset)$, this is the first condition in the criterion, and it is illustrated in Figure \ref{fig:example}(a).
 
We consider the set  $\dset$ (described below) of all monomials that are strictly $\pordering$  some element in $\cset$,
\begin{align}
	\nonumber	\dset &= \{ \beta \mid \exists \alpha \in \cset \text{ such that } \alpha \neq \beta \text{ and } \beta \pordering \alpha \} \\  &=\{\ell^j c^{m_1} \barc^{m_2}  \mid j \in \{0,1\},  j+ m_1 +m_2 \leq \EXP(n)  \}.
\end{align}

If  $\dset\ \cap\ \nonresidual(\fset)=\emptyset$ then we get $\cset\subseteq \minel_{\pordering} (\nonresidual(\fset)) $. 
For any monomial ordering $\monoorder$ and $\beta \in \dset$, if there is no polynomial $f \in \idealRes_{\fset}(\beta)$ such that $f \monostrict \beta$,  then $\dset \cap \nonresidual(\fset)=\emptyset$. Infact we show in Lemma \ref{lemma:TerminIdealres}, that if is sufficient to show that there is no monomial $\alpha \in \monoRes_{\fset}(\beta)$. This is the second condition of our criterion, this is illustrated in the Figure \ref{fig:example}(b). { In Theorem \ref{thm:main}, we state and prove  that the criterion described above,  ensures   $| \GB_{\monoorder}(\fset) | \geq \EXP(n)$}.

\begin{theorem}\label{thm:main}
	Let $\monoorder$ be any monomial ordering on $\monomials$, such that
	\begin{enumerate}[label=\roman*.,ref=\roman*]
		\item For every $\alpha \in \cset$, $s \monostrict \alpha$. \label{item:criteria:one}
			\item\label{item:criteria:two}For every $\beta \in \dset$, if $t \in \monoRes_{\fset}(\beta)$, then $ \beta \monoorder t$. 
	\end{enumerate}

Then, there are  at least $\EXP(n)$ many polynomials in $\GB_{\monoorder}(\fset)$ each with degree $\EXP(n)$.
 With this,  we also obtain  $| \GB_{\monoorder}(\fset) | \geq \EXP(n)$.

\end{theorem}

We postpone the proof of the theorem to the end of this section, we first establish $f\in \idealRes_\Poly(\beta)$ implies  $\termsof(f) \cap  \monoRes_\Poly(\beta) \neq \emptyset$. Hence, it is sufficient to consider  $\monoRes_{\fset}(\beta) $ instead of $\idealRes_{\fset}(\beta)$  in the above theorem.

\begin{lemma}\label{lemma:TerminIdealres}
	Let $\Poly\subseteq\F_2[X]$ be a set of binomials. For some $\alpha\in X^\oplus$ and $f\in \F_2[X]$, if $f\in \idealRes_\Poly(\alpha)$, then $ \termsof(f) \cap  \monoRes_\Poly(\alpha) \neq \emptyset$.
\end{lemma}

\begin{proof}
	For any monomial $\alpha \in X^\oplus$ and $f\in \idealRes_\Poly(\alpha)$,  suppose $\alpha\in \termsof(f)$, then clearly $ \alpha \in \termsof(f) \cap  \monoRes_\Poly(\alpha) \neq \emptyset$.
	
	Suppose $\alpha\notin \termsof(f)$, then we can write $\alpha+f= \sum_{i=1}^k\alpha_i \poly_i$, where $\alpha_1,\dots,\alpha_k \in X^{\oplus}$ and $\poly_1,\dots,\poly_k \in \Poly$. We need to prove that $t\in \monoRes_\Poly(\alpha)$ (equivalently $\alpha+t\in \Ideal{\Poly}$) for some $t\in \termsof(f)$. We do this by inducting on $k$.
	
	\paragraph{Base case}
	Suppose $k=1$, we have $\alpha+f=\alpha_1\poly_1$. Since $\poly_i$'s are binomials, we can write $\poly_i=\beta_{i1}+\beta_{i2}$ where $\beta_{i1}$ and $\beta_{i2}$ are monomials. Therefore, $\alpha+f=\alpha_1\beta_{11}+\alpha_1\beta_{12}$. Without loss of generality let $\alpha=\alpha_1\beta_{11}$, which implies $f=\alpha_1\beta_{12}$. Hence $f$ is a monomial, and  the result  follows.
	
	\paragraph{Case $k > 1$}
	We will assume  for any monomial $\alpha \in X^\oplus$ and $f\in \idealRes_\Poly(\alpha)$, if $\alpha+f=\sum_{i=1}^{k}\alpha_i \poly_i$ then there is a term $t\in \termsof(f)$ such that $t\in \monoRes_\Poly(\alpha)$ and we will prove the same for $k+1$.

	Assume $\alpha+f=\sum_{i=1}^{k+1}\alpha_i \poly_i$, for some monomial $\alpha$  and  $f\in \idealRes_\Poly(\alpha)$. Then we can write $\alpha=\alpha_j\beta_{j1}$ for some $j\in[1,k+1]$ (since $\alpha \notin \termsof(f)$). Then, $\alpha+f=\alpha_j\beta_{j1}+\alpha_j\beta_{j2}+\sum_{i=1,i\neq j}^{k+1}\alpha_i \poly_i$,   by rearranging we get,  $ \alpha_j\beta_{j2}+f = \sum_{i=1,i\neq j}^{k+1}\alpha_i \poly_i\in \Ideal{\Poly}$. Clearly, $f \in \idealRes_{\Poly}(\alpha_j\beta_{j2})$, by induction hypothesis we have $t \in \termsof(f)$ such that $ t\in \monoRes_\Poly(\alpha_j\beta_{j2})$ (that is, $\alpha_j\beta_{j2}+t\in \Ideal{\Poly}$). 
	
	Now note that  $\alpha_j\poly_{j} = \alpha_j\beta_{j1}+\alpha_j\beta_{j2} =\alpha+\alpha_j\beta_{j2}$ and hence $\alpha +t\in \Ideal{\Poly}$, this implies that $t \in  \monoRes_\Poly(\alpha)$ as required.
\end{proof} 

Let $\monoorder$ be any monomial ordering that satisfies the properties in Theorem \ref{thm:main}. Note that $\dset\cap \nonresidual(\fset)=\emptyset$ follows from the assumption (\ref{item:criteria:two}) of Theorem \ref{thm:main}. In order to show $\cset\subseteq \nonresidual(\fset)$, we prove that for any  $\monomial \in \cset$, $\monomial + s \in \Ideal{\fset}$. This is sufficient since  $s \monostrict \monomial$.   First we will establish that $b_{4n}^{\EXP(n)}c_{4n}f_n + c_{4n}s_n \in \Ideal{\pset}$ and $ {\barb_{4n}^{\EXP(n)}\barc_{4n}\barf_n + \barc_{4n} \bars_n\in \Ideal{\barP}}$ in the following proposition, the proof of this proposition can be found in  Section \ref{sec:missingproofs}. Subsequently in Proposition \ref{prop:subseteq}, we will show for any $\monomial \in \cset$, $\monomial + s \in \Ideal{\fset}$,  the proof is summarized in Figure \ref{prop:subseteq}.

\begin{proposition}\label{prop:mayr1}  With notations as above, we have the following:
\[
\begin{array}{lcr}
	{b_{4n}^{\EXP(n)}c_{4n}f_n + c_{4n}s_n \in \Ideal{\pset}} \qquad& \text{ and }&\qquad {\barb_{4n}^{\EXP(n)}\barc_{4n}\barf_n + \barc_{4n} \bars_n\in \Ideal{\barP}}.
\end{array}
\]

\end{proposition}

\begin{proposition}\label{prop:subseteq}
	For any $m_1 + m_2 = \EXP(n)$,  $\ell c^{m_1}\barc^{m_2}  + s \in  \Ideal{\fset}$.
\end{proposition}

\begin{proof}
	  By Proposition \ref{prop:mayr1} we have $ \barb_{4n}^{\EXP(n)}\barc_{4n}\barf_n + \barc_{4n}\bars_n  \in \Ideal{\fset}$, combining this with $ \barc_{4n}\bars_n + s \in \fset$,  we obtain 
	$$\p_1 = \barb_{4n}^{\EXP(n)} \barc_{4n}\barf_n + s \in \Ideal{\fset}.$$
	
	\noindent
	Since $\barc_{4n}\barf_n + c_{4n} s_n\in \fset$, we obtain $  \barb_{4n}^{\EXP(n)} c_{4n} s_n+ s \in \Ideal{\fset}$.
	 $$\p'_2= \barb_{4n}^{\EXP(n)-1}c_{4n}s_n b  + s\in \Ideal{\fset}.$$ 
	 Using this argument   $m_1$ times, we obtain 
	 $$\p_2= \barb_{4n}^{\EXP(n)-m_1}c_{4n}s_n b^{m_1}  + s\in \Ideal{\fset}.$$ 
	\noindent
	Similarly to the above reasoning, since  $ \barb_{4n}c_{4n}s_n + c_{4n}s_n\barb \in \fset$, we get 
	$$\p_3=  c_{4n} s_n  b^{m_1}\barb^{m_2}  + s\in \Ideal{\fset}.$$
	\noindent
	By Proposition \ref{prop:mayr1}, we know that  $ b_{4n}^{\EXP(n)} c_{4n} f_n + c_{4n} s_n\in \Ideal{\fset}$, hence we obtain 
	$$\p_4= b_{4n}^{\EXP(n)}c_{4n} f_n b^{m_1}\barb^{m_2} + s \in \Ideal{\fset}.$$
	\noindent
	Now since $\ c_{4n} f_n+ \ell \in \fset$ , we obtain 
	
	$$\p_5=  b_{4n}^{\EXP(n)}  b^{m_1}\barb^{m_2}\ell  + s \in \Ideal{\fset}.$$
	Finally since $ b_{4n} b \ell + \ell c,\  \barb b_{4n} \ell+ \ell\barc \in \fset$,  we have $\ell  c^{m_1}\barc^{m_2}+ s \in \Ideal{\fset}$.
\end{proof}
 
 \noindent

 \begin{figure}[H]
 	\includegraphics[scale=0.32]{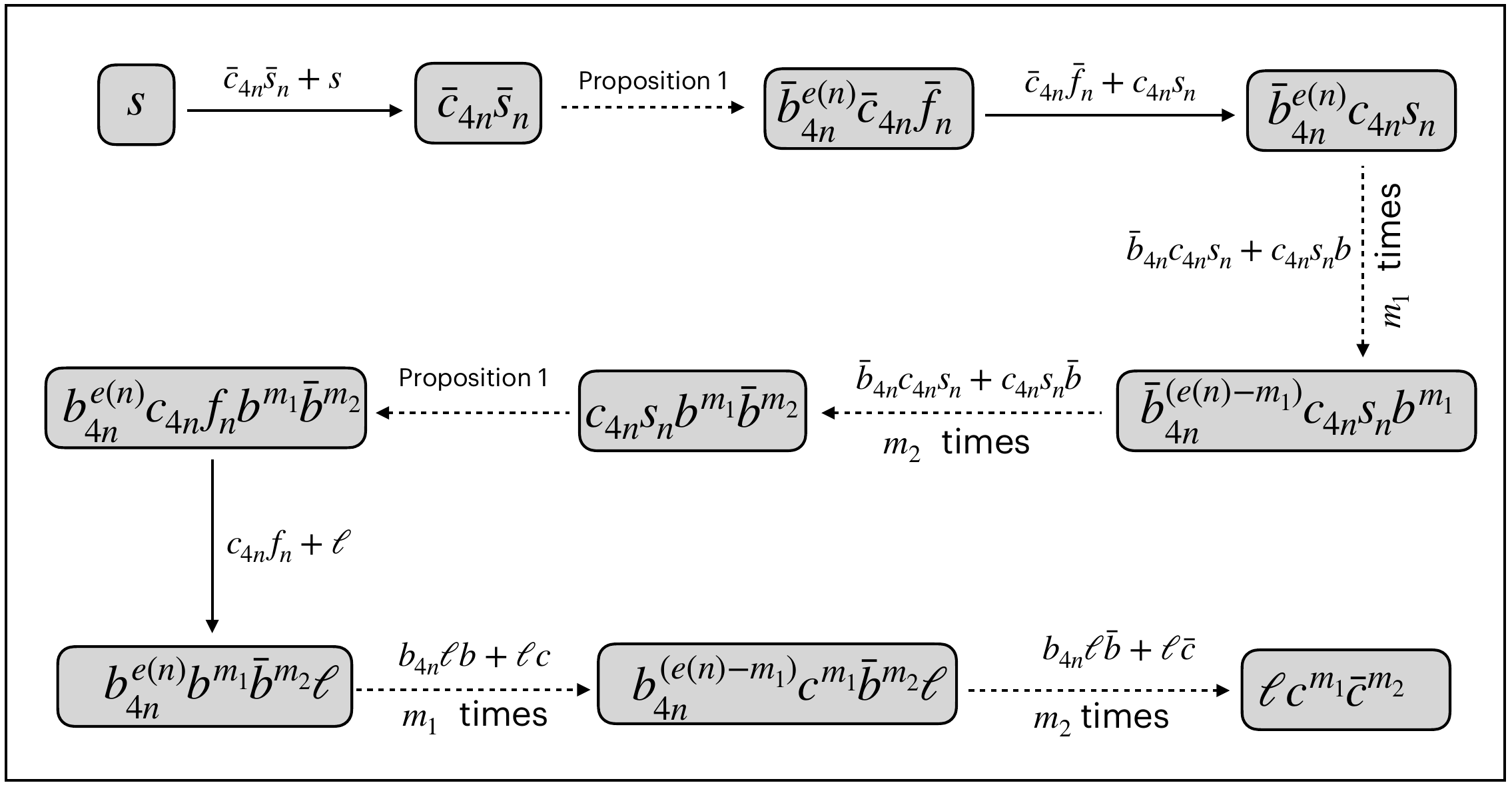}
 	\caption{ $s + \ell c^{m_1} \barc^{m_2} \in \Ideal{\fset}$}
 \end{figure}

\begin{lemma}\label{lem:csetinsf}
	For any monomial ordering $\monoorder$ and for every $\alpha \in \cset$, if $s \monostrict \alpha$ then $\cset \subseteq  \nonresidual(\fset)$.
\end{lemma}	
	
\begin{proof}
By Proposition \ref{prop:subseteq}, for any $\alpha \in \cset,\ \alpha+s\in \Ideal{\fset} $. That is, $\alpha\xrightarrow[\fset]{} s$. Since $s\monostrict \alpha$, we have $\alpha \in \nonresidual(\fset)$ and hence $\cset \subseteq \nonresidual(\fset)$.
\end{proof}

\begin{proof}[Proof of Theorem \ref{thm:main}]
By Lemma \ref{lem:csetinsf}, $\cset \subseteq \nonresidual(\fset)$. Also we have, for $\alpha\in \cset$, if $\beta\pordering \alpha$ then $\beta \in \dset$. 
 If $\ffun \in \idealRes_{\fset}(\beta)$, then there is a term $t \in \termsof(f)$ such that $ t \in \monoRes_{\fset}(\beta)$ by Lemma \ref{lemma:TerminIdealres}. Notice that if $ \beta \monoorder t$ then  $\beta \monoorder \ffun$, hence under the assumption in (\ref{item:criteria:two}), we also get $\beta   \monoorder \ffun$.
So, we obtain $\cset \subseteq \minel_\pordering(\nonresidual(\fset))$, and hence  $ |\minel_\pordering(\nonresidual(\fset))|\geq |\cset| \geq \EXP(n)$. 

Since all elements in $\cset$ have degree greater than $\EXP(n)$, and since $\minel_{\pordering}(\nonresidual(\fset))=\htermsof(\GB(\fset))$ (by Theorem \ref{thm:cardGB}), for each $\alpha\in \cset$, there exists a polynomial $p\in\GB(\fset)$ such that $\htermsof(p)=\alpha$. Hence, there are  at least $\EXP(n)$ many polynomials in $\GB(\fset)$ each with degree $\EXP(n)$.
\end{proof}

\section{Large cardinality for reduced Gr\"obner basis } \label{sec:MonoOrder}

In this section, we use the notations from Section \ref{sec:SEGB}. Now, we investigate various monomial orderings that satisfy the criterion of Theorem \ref{thm:main}. We will first prove Proposition \ref{prop:snnotdivD}, which reasons about certain interesting properties about the ideal residues of polynomials in $\dset$. For proving Proposition \ref{prop:snnotdivD}, we state and prove certain lemmas  that will be useful below.  By abuse of notation, we say that a binomial $f \in \fset$ is \emph{applicable} from a monomial $\alpha$ if there is a $t  \in \termsof(f)$ such that $\divides{t}{\alpha}$. Note that if $f $ is applicable from a monomial $\alpha$, then there is an $\alpha'$ such that $\alpha+ \alpha' = m f$ for some monomial $m$, in this case we say $\alpha$ \emph{reduces} to $\alpha'$, we also use this  to sometimes refer to reduction due to multiple steps.

\begin{example}
Let $\alpha=\ell c^{2} \barc^{2}$,  consider the binomial $b_{4n}\ell b+\ell c\in \fset$. This binomial is applicable from $\alpha$, since $\divides {\ell c}{\alpha}$. Since $\alpha+\ell c \barc^2 b_{4n} b= c\barc^2(b_{4n}\ell b+\ell c)$, we say $\alpha$ reduces to $\ell c \barc^2 b_{4n} b$. 

\end{example}

\begin{lemma}\label{lem:seq}
	Let   $\Poly \subseteq \Z_2[X]$ be a set of binomials and  $\monomial,\monomial' \in X^\oplus$ be any two monomials such that $\monomial + \monomial' \in \Ideal{\Poly}$. Then we can find a sequence of monomials $\monomial_0, \ldots , \monomial_r \in X^\oplus$ such that
	
	\[
	\monomial = \monomial_0,\ \monomial'=\monomial_r,\text{ and } \forall i \in [0 ,r -1],\ \monomial_i + \monomial_{i+1} =  \beta_i \poly_i
	\]

	for some $\beta_i \in X^\oplus$ and $\poly_i \in \Poly$.
\end{lemma}
\begin{proof}
	Since $\monomial + \monomial' \in \Ideal{\Poly}$, we have  $\monomial + \monomial' = \Sigma_{j=1}^m a_j  \poly_j$ for some $a_1, \ldots, a_m \in X^\oplus$, $\poly_1, \ldots, \poly_m \in \Poly$.  We prove the statement by inducting on $m$.
	
	Suppose $m = 1$, then we have nothing to prove. Assume that $\monomial \neq \monomial'$, then  $ \monomial \in \termsof(a_{l} \poly_{l})$ for some $l \in [1,m]$. Since $\poly_{l} = \bar{\monomial}_1 + \bar{\monomial}_2$ for some $ \bar{\monomial}_1 , \bar{\monomial}_2 \in X^\oplus$, then, without loss of generality, let $\alpha=a_{l}\bar{\monomial}_1$ and let $ a_{l}\bar{\monomial}_2= \monomial''$.
	Then $\monomial''+\monomial' =  \Sigma_{j=1, j \neq l}^m a_j \poly_j $. Now by the induction hypothesis there is a monomial sequence for $\monomial''+\monomial'\in \Ideal{\Poly}$. Combining with $\monomial$, we get the required result.
\end{proof}

\begin{lemma}\label{lem:mayrmeyer}
	For any $ j \in [1,4]$ and $\alpha \in   \monoRes_\pset(c_{jn} s_n)$,  if $\divides{f_n} {\alpha}$ then $\alpha = b_{jn}^{\EXP(n)}c_{jn}f_n $.
\end{lemma}

\noindent Proof of Lemma \ref{lem:mayrmeyer} can be found in Section \ref{sec:missingproofs}.

\begin{lemma}\label{lem:preserve}
	Let $\beta \in \dset$ be such that  $\beta = \ell c^{m_1} \barc^{m_2}   $.  Let $\gamma \in \monoRes_{\fset}(\beta)$ and $\divides{\ell}{\gamma}$, then  $\gamma = \ell c^{m'_1} \barc^{m'_2} b^{m''_1} \barb^{m''_2} {b_{4n}^{m''_1 + m''_2}}$ where $m_1 = m'_1+ m''_1$ and $m_2 = m'_2 + m''_2$. 
\end{lemma}

\begin{proof}
	
This follows from a simple observation of the polynomials. The only polynomials that involves $c,\ b$ and $\barc,\ \barb$ are   $	b_{4n}\ell b+\ell c$ and $\   b_{4n}\ell \barb + \ell \barc $. Moreover the only other polynomial involving the variable $\ell$ is $c_{4n} f_n+ \ell$. In each of these polynomials the sum of the degrees of $c,b$ and $\barc,\barb$ remains invariant.	Moreover, application of these polynomials ensures that the degree of $b_{4n}$ is equal to the sum of the  degrees of $b$ and $ \barb$ in any $\gamma \in \monoRes_{\fset}(\beta)$.
\end{proof}

\begin{lemma}\label{lem:notdivsn}
 Let $\beta = c_{4n} f_n {b_{4n}^{m}} \beta' $ with $  m < \EXP(n)$ and $\beta' \in \{c, \barc, b, \barb \}^\oplus$, then for any $\gamma \in \monoRes_{\fset \setminus \sset_\ell}(\beta)$,  $\notdivides{ c_{4n}s_n}{\gamma}$, where $\sset_\ell = \{ b_{4n}\ell b+\ell c,  b_{4n}\ell \barb + \ell \barc, c_{4n} f_n+ \ell\}$
\end{lemma}
\begin{proof}
	
	 Let us assume that $\divides{ c_{4n}s_n}{\gamma}$ for some $\gamma \in \monoRes_{\fset \setminus \sset_\ell}( \beta ) $.   Then by Lemma \ref{lem:seq}, there are $\alpha_0, \alpha_1, \dots, \alpha_m$ such that $\alpha_0 = \beta$, $\alpha_m = \gamma$ and $\forall i \in [0 ,m -1],\ \monomial_i + \monomial_{i+1} =  \beta_i \poly_i$ for some $\poly_i \in \fset \setminus \sset_\ell$.  We will assume that this sequence is minimal, that is for any $j < m$, $\notdivides{c_{4n}s_n}{\alpha_j}$.
	Notice that the polynomials from 
	$\sset = \{ \barb_{4n}c_{4n} s_n +c_{4n}s_n\barb, \barb_{4n}c_{4n}s_n +c_{4n}s_n b\}$ are guarded by $c_{4n}s_n$ 
	and hence, for any $j \in [0,m-1]$,  $\poly_j \not\in \sset$.   
	Moreover, there is no $j \in [0,m-1]$ such that  $\poly_{j} = \barc_{4n}\barf_n +c_{4n}s_n$, since this is the only polynomial that can reduce monomials from $V^\oplus$ to monomials in $\bar{V}^\oplus$. A similar reasoning can be used to show that  $\poly_j \neq \barc_{4n} \bars_n+ s$.
	
	 Hence, for any $j \in [0,m-1]$, $ \poly_j \in \pset$ and  $\gamma \in \monoRes_{\pset}(\beta)$ such that $\divides{ c_{4n}s_n}{\gamma}$. 
	 Since $m < \EXP(n)$,  this contradicts Lemma \ref{lem:mayrmeyer}.
\end{proof}

\begin{lemma}\label{lem:notdivbarsn}
	Let $\beta$ be a monomial such that $\beta = b^{m_1} \barb^{m_2} c_{4n} f_n   {b_{4n}^{m}}\beta' $ for some $\beta' \in \{c, \barc \}^\oplus$ and $m = m_1 + m_2$. 	
	For any monomial $\gamma \in \monoRes_{\fset \setminus \sset_\ell}(\beta)$, if  $m < \EXP(n)$, then $\notdivides{ \barc_{4n}\bars_n}{\gamma}$.
\end{lemma}
\begin{proof}
			Let us assume that there is a $\gamma \in \monoRes_{\fset \setminus \sset_\ell}(\beta)$ such that $m < \EXP(n)$, and  $\divides{ \barc_{4n}\bars_n}{\gamma}$. Then by Lemma \ref{lem:seq}, there are $\alpha_0, \alpha_1, \dots, \alpha_m$ such that $\alpha_0 = \beta$, $\alpha_m = \gamma$ and $\forall i \in [0 ,m -1],\ \monomial_i + \monomial_{i+1} =  \beta_i \poly_i$ for some $\poly_j \in {\fset \setminus \sset_\ell}$, we assume  that the sequence is minimal.
			
			Our strategy is to show that if $\divides{\barc_{4n}\bars_n}{\gamma} $, then there is an $i < m$ such that $\divides{c_{4n}s_n}{\alpha_i}$. From this, we will obtain a contradiction to Lemma \ref{lem:notdivsn}. We will first show that for any $j \in [0,m-1]$, $\poly_j \neq \barc_{4n} \bars_n+ s$.  Notice that this is the only polynomial involving the term $s$. Hence if this polynomial was used in the reduction, then we would also obtain a shorter sequence,  which would contradict the minimality of the assumed sequence.
			
			Notice that the only polynomial in $\gset \cup \pset$ such that one of its terms is divisible by $\barc_{4n} \bars_n$ is $\barc_{4n} \bars_n+ s$, as argued above, this polynomial does not occur in the sequence. 
			As in the proof of the previous lemma, a polynomial from $\barP$ can be used for the reduction only after using the polynomial $ \barc_{4n}\barf_n +c_{4n}s_n$, and from this we will obtain an $\alpha_i$ for some $i\in [0,m-1]$ such that $\divides{\alpha_i}{c_{4n}s_n}$.   
\end{proof}

\begin{proposition}\label{prop:snnotdivD}
	Let $\alpha  \in \dset$, and  $\gamma \in \monoRes_{\fset}( \alpha )$, then the following holds.
	\begin{enumerate}
		\item If $\notdivides{\ell}{\alpha} $, then $\monoRes_{\fset}( \alpha ) = \{ \alpha\}$  \label{prop:item:one}
		\item $\notdivides{ c_{4n}s_n}{\gamma}$ and $\notdivides{ \barc_{4n}\bars_n}{\gamma}$. \label{prop:item:two}
		\item $\monoRes_{\fset}( \alpha )\cap \set{\ell c^{r_1}\barc^{r_2}\suchthat r_1,r_2\in \N}=\set{\alpha}$.
		 \label{prop:item:three}
 		\item If $\gamma \neq \alpha$, $\deg(\gamma) > \deg(\alpha)$. \label{prop:item:four}
	\end{enumerate}
\end{proposition}

\begin{proof}

\textsf{Proof of \ref{prop:item:one}:} Recall that $ \dset = \{\ell^j c^{m_1} \barc^{m_2}  \mid j \in \{0,1\},  j+ m_1 +m_2 \leq \EXP(n) \} $. If $\notdivides{\ell}{\alpha}  $ then $\alpha = c^{m_1} \barc^{m_2} $ for some $m_1, m_2$ such that $m_1 + m_2 \leq \EXP(n)$. 
Since there is no
polynomial  $f \in \fset$ that is applicable from $\alpha$,   there is no $\alpha' \in \varset^\oplus$ such that $\alpha + \alpha' \in \Ideal{\fset}$. From this, the proof is immediate.

\medskip

\textsf{Proof of \ref{prop:item:two}:} If $\alpha=c^{m_1} \barc^{m_2}$, then  by item \ref{prop:item:one} of this lemma, we have $\gamma=\alpha $ and hence $c_{4n} s_n$ and  $\barc_{4n}\bars_n$ do not divide $\gamma$.

We assume that $\alpha=\ell c^{m_1}\barc^{m_2}$, then we have $m_1 + m_2 < \EXP(n) $. For the sake of contradiction, let us assume that $\divides{ c_{4n}s_n}{\gamma}$ for some $\gamma \in \monoRes_{\fset}( \alpha ) $.   Then by Lemma \ref{lem:seq}, there are $\alpha_0, \alpha_1, \dots, \alpha_m$ such that $\alpha_0 = \alpha$, $\alpha_m = \gamma$ and $\forall i \in [0 ,m -1],\ \monomial_i + \monomial_{i+1} =  \beta_i \poly_i$ for some $\poly_i \in \fset$. 
Since $\divides{\ell}{\alpha}$, we can find an $\alpha_k$ with $k$ being maximum such that $\divides{\ell}{\alpha_k}$.  Then by Lemma \ref{lem:preserve} we have $\alpha_k =  \ell c^{m'_1} \barc^{m'_2} b^{m''_1} \barb^{m''_2} {b_{4n}^{m''_1 + m''_2}}$ where $m_1 = m'_1 + m''_1$ and $m_2 = m'_2 + m''_2$. 
Clearly $\poly_{k} = \ell+f_nc_{4n}$ and hence $\alpha_{k+1} =  f_n c_{4n} c^{m'_1} \barc^{m'_2} b^{m''_1} \barb^{m''_2} {b_{4n}^{m''_1 + m''_2}}$. Then we also have $\gamma \in \monoRes_{\fset \setminus \sset_\ell}(\alpha_{k+1})$, this contradicts Lemma \ref{lem:notdivsn}. The second part follows directly from Lemma \ref{lem:notdivbarsn}.

\medskip
\textsf{Proof of \ref{prop:item:three}:}  This follows from a simple observation of the polynomials. The only polynomials that involve $c$ and $\barc$ are from  $	\{ b_{4n}\ell b+\ell c,\   b_{4n}\ell \barb + \ell \barc\} $. In each of these polynomials the sum of the degrees of $c,b$ and $\barc,\barb$ remains invariant. We provide a formal proof.    Let 
\[ 
	\sset_1 = \monoRes_{\sset_\ell \setminus \{f_nc_{4n} + \ell \} }(\alpha)\qquad  \quad 
 \sset_2 =   \monoRes_{\sset_\ell  }(\alpha) \setminus \sset_1
 \]

 Observe that by Lemma \ref{lem:preserve}, for any $\beta\in S_1$,  $\beta$ is of the form $\ell c^{m'_1} \barc^{m'_2} b^{m''_1} \barb^{m''_2} {b_{4n}^{m''_1 + m''_2}}$ where $m_1 = m'_1+ m''_1$ and $m_2 = m'_2 + m''_2$. 
  Note that $\monoRes_{\fset}( \alpha )=\sset_1 \cup \sset_2 \ \bigcup_{\beta\in \sset_2} \monoRes_{\fset \setminus \sset_\ell}( \beta )$. Let $\alpha= \ell c^{m_1}\barc^{m_2}$ and let $\gamma\in \monoRes_{\fset}( \alpha )\cap \set{\ell c^{r_1}\barc^{r_2}\suchthat r_1,r_2\in \N}$.  Since $\divides{\ell}{\gamma}, $ we can observe that $\gamma \in \sset_1$ and from the definition of $S_1$, we get $m_1''=m_2''=0$. This gives $r_1=m_1$ and $r_2=m_2$, and hence $\gamma=\alpha$.
 
\medskip
\textsf{Proof of \ref{prop:item:four}:} Let $\alpha=\ell c^{m_1}\barc^{m_2}$ (which is the only possibility by \ref{prop:item:one}). Now, consider the set $\sset_3 =  \monoRes_{\fset}( \alpha ) \setminus \sset_1 \cup \sset_2$. We will separately argue that the degree is large in each of the sets $\sset_1, \sset_2,\sset_3$.

Let   $\beta \in \sset_1$, then $\beta = \ell c^{m'_1} \barc^{m'_2} b^{m''_1} \barb^{m''_2} {b_{4n}^{m''_1 + m''_2}}$, for some $m''_1, m''_2$ such that  $m''_1 + m''_2 > 0$. Further, we have $m_1 +m_2 = m'_1 + m'_2+ m''_1 + m''_2$, which implies $\deg(\beta) > \deg(\alpha)$.

For any $\beta \in \sset_2$, it is obtained by replacing $\ell$ in some $\beta' \in \sset_1$ by $f_n c_{4n}$, and hence we have $\deg(\beta) > \deg(\alpha)$.

For any $\beta \in \sset_3$,
we have $\beta \in \monoRes_{\pset}(\beta')$ for some $\beta' \in \sset_2$, this is since no polynomials from $\gset \cup \barP$ can be used in the reduction. Let $\beta' = f_n c_{4n} c^{m'_1} \barc^{m'_2} b^{m''_1} \barb^{m''_2} {b_{4n}^{m''_1 + m''_2}}$, then  $\divides{ c_{4n}c^{m'_1} \barc^{m'_2} b^{m''_1} \barb^{m''_2}}{\beta}$. We only have to prove that $\beta \neq { c_{4n}c^{m'_1} \barc^{m'_2} b^{m''_1} \barb^{m''_2}}$. But this is easy to see since $f_n c_{4n} {b_{4n}^{m''_1 + m''_2}}$ never reduces to $c_{4n}$ through $\pset$.
\end{proof}

Now, we will prove that the following monomial orderings satisfy the criterion of Theorem \ref{thm:main}, as a result the corresponding reduced Gr\"obner basis for $\Ideal{\fset}$ is double exponential in size.

\subsection{Lexicographic ordering}\label{subsec:lex}
\begin{definition}\label{def:lex}
The \emph{lexicographic  ordering} ($\monoorder_{lex}$)  between monomials is defined as follows: Let $\alpha=x_1^{\degree_1}\cdots x_n^{\degree_n}$ and $\beta=x_1^{\degree'_1}\cdots x_n^{\degree'_n}$, then $\alpha \monoorder_{lex} \beta$ if and only if $\alpha = \beta$ or there is an $i$ with $1\leq i\leq n$ such that for all $1 \leq j < i$, $\degree_j = \degree'_j$ and $\degree_i < \degree'_i$.
\end{definition}

We will prove that the lexicographic ordering satisfies the criterion in Theorem \ref{thm:main}. 
Let us order the variables in $\varset$ lexicographically as follows. 
We recall   the set of all variables  below:
\begin{equation*}
	\begin{split}
		\varset  =\  &  V \cup  \bar{V}    \cup \{ s,\ell,c,\barc,b,\barb\}, \text{ where }  
		V =  \bigcup_{i=0}^n V_i,  \quad   \bar{V} =  \bigcup_{i=0}^n \bar{V}_i  \text{ and }\\
		& V_i = \{s_i,f_i,q_{1i},q_{2i},q_{3i},q_{4i},c_{1i},c_{2i},c_{3i},  c_{4i},b_{1i},b_{2i},b_{3i},b_{4i}  \}	\\
		& \bar{V}_i = \{\bars_i,\barf_i,\barq_{1i},\barq_{2i},\barq_{3i},\barq_{4i},\barc_{1i},\barc_{2i},\barc_{3i},  \barc_{4i},\barb_{1i},\barb_{2i},\barb_{3i},\barb_{4i}  \}	.\\
	\end{split}
\end{equation*} 

Let $S = \{ s_0, \ldots , s_n\}$ and $F = \{ f_0, \ldots, f_n\}$. For $j \in [1,4]$, let $Q_j = \{ q_{j0}, \ldots, q_{jn}\}$, $C_j = \{ c_{j0}, \ldots, c_{jn}\}$, $B_j= \{ b_{j0}, \ldots,  b_{jn}\}$, and $\bar{S}, \bar{F}, \barQ_j, \barC_j, \barB_j$ defined analogously by including an appropriate set of variables from $\bar{V}$. 

Between any  $V_i$ and $V_j$, where $0 \leq i < j \leq n$, we establish the following ordering.
We let   
$ s_i \monoorder_{lex}s_j, f_i \monoorder_{lex}f_j$,  $b_{ki}\monoorder_{lex} b_{kj},  c_{ki}\monoorder_{lex} c_{kj}, \text{ and } q_{ki}\monoorder_{lex} q_{kj}$ where $k \in [1,4]$. We establish a similar ordering between  $\bar{V_i}$ and $\bar{V_j}$. 
The overall  ordering is established as follows.

\begin{equation}\label{eq:lexorder}
	\begin{split}
		\{s\} & \monoorder_{lex} \{c\}  \monoorder_{lex} \{\barc\} \monoorder_{lex}  \{ \ell \}  \monoorder_{lex} \{b\} \monoorder_{lex} \{ \barb\} \monoorder_{lex}  S \monoorder_{lex}  F \monoorder_{lex}\\
		& C_1 \monoorder_{lex} 
		 \cdots \monoorder_{lex} C_4 \monoorder_{lex} B_1 \monoorder_{lex} \cdots \monoorder_{lex} B_4 \monoorder_{lex}  Q_1 \monoorder_{lex} \cdots \monoorder_{lex} \\ & Q_4
		\monoorder_{lex} \barS \monoorder_{lex} \bar{F} \monoorder_{lex} \barC_1 \monoorder_{lex} \cdots \monoorder_{lex} \barC_4 \monoorder_{lex}  \barB_1 \monoorder_{lex} \cdots \monoorder_{lex} \\ & \barB_4  \monoorder_{lex}  \barQ_1  \cdots \monoorder_{lex} \barQ_4.
	\end{split}
\end{equation}

Since $s$ is lexicographically less than all other variables, this ordering satisfies the first condition in Theorem \ref{thm:main}, and the second condition will be proved in the following proposition. 

\begin{proposition}\label{prop:lex}
For every $\beta \in \dset$, if $ \gamma \in \monoRes_{\fset}(\beta) $ then $ \beta \monoorder_{lex}  \gamma$.
\end{proposition}
\begin{proof}
Since $\beta \in \dset$, $\beta= \ell^j c^{m_1} \barc^{m_2}  $ where $j \in \{0,1\},$ with $  j+ m_1 +m_2 \leq \EXP(n)$. 
If $\notdivides {\ell}{\beta}$, then by Proposition \ref{prop:snnotdivD}-\ref{prop:item:one}, $\monoRes_{\fset}( \beta ) = \{ \beta\}$, and hence $ \gamma = \beta $.

If $\divides{\ell}{\beta}$ then assume for the sake of contradiction that $\gamma \monostrict_{lex} \beta$. Then from the lexicographic ordering of variables as defined
in \eqref{eq:lexorder}, it follows that no variables other than $s,c,\barc,\ell $ divide $\gamma$. We consider two cases.

Suppose $\notdivides{s}{\gamma} $. Then $\gamma $ is of the form $\ell c^{r_1} \barc^{r_2}$ where $r_2< m_2$ or ($r_2 = m_2$ and $r_1 < m_1$). But by Proposition \ref{prop:snnotdivD}-\ref{prop:item:three}, if $\gamma=\ell c^{r_1} \barc^{r_2}\neq \beta$, then $\gamma\notin \monoRes_{\fset}( \beta )$. So, this case is not possible.

Suppose $\divides{s}{\gamma} $, since  $\gamma \in \monoRes_{\fset}( \beta )$, by Lemma \ref{lem:seq} we have a monomial sequence $\alpha_0=\beta,\alpha_1,\ldots,\alpha_r=\gamma$ with $\alpha_i+\alpha_{i+1}=\beta_i\poly_i$, for all $i\in[0,r-1]$, and $\poly_i \in \fset$. Since $\barc_{4n} \bars_n+s\in \gset$ is the only polynomial in $\fset$ which contains $s$, we can conclude that $\divides{\barc_{4n} \bars_n}{\alpha_{r-1}} $.  
This contradicts  Proposition \ref{prop:snnotdivD}-\ref{prop:item:two}. 
\end{proof}

\begin{theorem}

	Let $\monolex$ be the monomial ordering defined as above. Then, there are at least $\EXP(n)$ many polynomials in $\GB_{\monolex}(\fset)$ each with degree $\EXP(n)$.
	With this,  we also obtain  $| \GB_{\monolex}(\fset) | \geq \EXP(n)$.

\end{theorem}

\subsection{ Degree lexicographic ordering}\label{subsec:deg}

\begin{definition}\label{def:deglex}
	The \emph{degree lexicographic ordering} ($\monodeg$)  between monomials is defined as follows: Let $\alpha,\beta\in \monomials$, then $\alpha\monodeg \beta$ if and only if $\deg(\alpha)< \deg(\beta)$ or $(\deg(\alpha)=\deg(\beta)$ and $\alpha\monolex \beta)$, where $\monolex $ is as in Subsection \ref{subsec:lex}. 
\end{definition}

In this case, it is obvious that  $s$ is degree lexicographically less than all elements in $\cset$. Also, by Proposition \ref{prop:snnotdivD}-\ref{prop:item:four}, for every $\beta \in \dset$, if a monomial $\gamma \in \idealRes_{\fset}(\beta)$ and $\gamma \neq \beta$,  then $\deg(\gamma) > \deg(\beta) $. Hence, we obtain the following  theorem.

\begin{theorem}
	
	Let $\monodeg$ be the monomial ordering defined as above. Then, there are  at least $\EXP(n)$ many polynomials in $\GB_{\monodeg}(\fset)$ each with degree $\EXP(n)$.
	With this,  we also obtain  $| \GB_{\monodeg}(\fset) | \geq \EXP(n)$.
	
\end{theorem}

\begin{remark}
	    The case of reverse degree lexicographic ordering is similar to that of the degree lexicographic ordering.
\end{remark}

\subsection{Weighted ordering}\label{subsec:wt}

We consider a special case of weighted ordering with the definition as presented in \cite[Exercise.11,p.74]{CLO1997} and show that the reduced Gr\"obner basis in this case is also double exponential in size. 

\begin{definition} 
	Given a set of variables $\vset$,  we say $\wtmap:\vset \to \R_{> 0}$ is a \emph{weight map} if it is injective and the set $ \{\wtmap(x) \mid x \in \vset\}$ is  linearly independent over $\Q$. Given a weight map $\wtmap$ and a monomial $\gamma= \prod_{x\in\vset} x^{\degree_{x}}$, we let $\wtdeg_\wtmap(\gamma) =\sum \limits_{x \in \vset} \degree_x\wtmap(x)$ and refer to it as the weighted degree of $\gamma$ with respect to $\wtmap$. For $\alpha,\beta\in \vset^{\oplus}$, we say $\alpha\monoorder_\wtmap \beta $ if and only if $\alpha=\beta$ or $\wtdeg_\wtmap(\alpha)< \wtdeg_\wtmap(\beta)$. We denote by $ \monoorder_\wtmap$ the \emph{weighted ordering} determined by $\wtmap$. 
	
\end{definition}

Recall that $\varset$ is the set of all variables that appear in $\fset$, we wish to find an appropriate weight map for it.  Let $\wtmap:\varset \to \R_{> 0}$ be a weight map such that it satisfies the following conditions.
\begin{enumerate}
	\item $\wtmap(s) < \wtmap(c)< \wtmap(\barc) < \wtmap(\ell) < \wtmap(b) < \wtmap(\barb) $.
	\item For any $x \in \varset \setminus \{s,c,\barc,\ell,b,\barb \} $ and $y \in  \{s,c,\barc,\ell,b,\barb \} $, $  \wtmap(y)< \wtmap(x)$.
\end{enumerate}

We claim that  the weighted ordering $ \monoorder_\wtmap$ satisfies the criterion specified in Theorem \ref{thm:main}. Since $\wtmap(s)$ is the least weight among all variables, it is easy to see that $s\monoorder_{\wtmap} \alpha$ for any $\alpha \in \cset$. This ensures the first condition. We next show that the second condition is also satisfied by $\monoorder_{\wtmap} $.

\begin{lemma}\label{lem:wtdeg}
For every  $\beta \in \dset$, if $\gamma\in \monoRes_{\fset}(\beta ) $ then $\beta\monoorder_{\wtmap}\gamma$.
\end{lemma}
\begin{proof}
Let $\beta \in \dset$, then $\beta= \ell^j c^{m_1} \barc^{m_2}  $, where $j\in \set{0,1}$ and $j+m_1 +m_2\leq \EXP(n)$. If $\notdivides {\ell}{\beta}$, then by Proposition \ref{prop:snnotdivD}-\ref{prop:item:one}, $\monoRes_{\fset}( \beta )= \{ \beta\}$,  hence $\gamma = \beta$ and $\beta \monoorder_{\wtmap}\gamma$.

If $\divides{\ell}{\beta}$, then consider the sets $\sset_1,\sset_2$ and $\sset_3$ as in the proof of Proposition \ref{prop:snnotdivD}. Notice that $\monoRes_{\fset}(\beta )  =S_1\cup S_2 \cup S_3$.

Suppose $\gamma \in \sset_1 $, then $\gamma = \ell  c^{m'_1} \barc^{m'_2}  b^{m''_1} \barb^{m''_2}{b_{4n}^{m''_1 + m''_2}}$ such that
 $m_1= m'_1 +  m''_1$ and $m_2= m'_2+m''_2$  and $m''_1 + m''_2 > 0$. 
 Then $\wtdeg_\wtmap(\beta)= \wtmap(\ell) +m_1 \wtmap(c)+m_2 \wtmap(\barc)$ and $\wtdeg_\wtmap(\gamma)=\wtmap(\ell)+m'_1 \wtmap(c)+m'_2 \wtmap(\barc)+m''_1 \wtmap(b)+m''_2 \wtmap(\barb)+(m''_1+m''_2) \wtmap(b_{4n})$. Since  $\wtmap(c)< \wtmap(\barc) < \wtmap(b) < \wtmap(\barb) $, it is easy to see that $ \wtdeg_\wtmap(\beta)<\wtdeg_\wtmap(\gamma)$.

Suppose $\gamma \in \sset_2 $, then clearly $\gamma= c_{4n}f_nc^{m'_1} \barc^{m'_2} b^{m''_1} \barb^{m''_2}$. Since $\wtmap(\ell) < \wtmap(c_{4n})$, we will obtain $\beta \monoorder_{\wtmap} \gamma$ and we are done. 
 The case when $\gamma \in \sset_3 $ is similar.
\end{proof}

\begin{theorem}
	
	Let $\monoorder_{\wtmap}$ be the monomial ordering defined as above. Then, there are at least $\EXP(n)$ many polynomials in $\GB_{\monoorder_{\wtmap}}(\fset)$ each with degree $\EXP(n)$.
	With this, we also obtain  $| \GB_{\monoorder_{\wtmap}}(\fset) | \geq \EXP(n)$.
	
\end{theorem}

\section{Missing Proofs}\label{sec:missingproofs}

\subsection{Proof of Proposition \ref{prop:mayr1}}

The proof we present here is based on the proof in \cite[sec.6]{MM1982}. 
\begin{proof}
We will use induction to prove the proposition, in fact we will only prove the first part, as the second part is symmetric. 
	That is, we only show ${b_{4n}^{\EXP(n)}c_{4n}f_n + c_{4n}s_n \in \Ideal{\pset}}$.
	Towards this, we prove the stronger statement that for any $i \in [1,4]$ and $m \leq n$, ${b_{im}^{\EXP(m)}c_{im}f_m + c_{im}s_m \in \Ideal{\pset}}$. We prove this by inducting on $m$.
	We refer readers to Figure \ref{fig:meyermeyer} for the summary of the proof.

	\paragraph{ Case  $m = 0$}This follows directly  from the set of polynomials in $\pset_0 $.
	
	\paragraph{ Case  $m +1$} We fix an $i \in [1,4]$ and show that $b_{i{(m+1)}}^{\EXP(m+1)}c_{i{(m+1)}}f_{(m+1)} + c_{i{(m+1)}}s_{(m+1)} \in \Ideal{\pset}$.  We assume by the induction hypothesis that for each $i  \in [1,4]$, {$b_{i{m}}^{\EXP(m)}c_{i{m}}f_{m} + c_{i{m}}s_{m} \in \Ideal{\pset}$}. 
	\vskip 15pt
	
	\noindent
	 Consider the polynomial  $ q_{1{(m+1)}}c_{1{m}}s_{m}+s_{m+1} \in \pset$. By multiplying $c_{i(m+1)}$, we have the following
	 {\[  \p_1 =      q_{1{(m+1)}}c_{i{(m+1)}}c_{1{m}} s_{m} + c_{i{(m+1)}}s_{(m+1)}  \in  \Ideal{\pset}	\]}
	
	\noindent
	By induction hypothesis, we have {$ b_{1{m}}^{\EXP(m)}c_{1{m}} f_{m}+c_{1{m}}s_{m} \in \Ideal{\pset}$} and hence we have 
	{\[ \p_2 = 	 q_{1{(m+1)}}c_{i{(m+1)}}b_{1{m}}^{\EXP(m)} c_{1{m}} f_{m}+ c_{i{(m+1)}}s_{(m+1)} \in \Ideal{\pset} \]}
	
	\noindent	
	Since $q_{2(m+1)}c_{2m}s_{m}+q_{1(m+1)}b_{1m}c_{1m}f_{m} \in \pset$, we have 
 {\[	\p_3 = q_{2(m+1)} c_{i{(m+1)}} b_{1{m}}^{(\EXP(m)-1)} c_{2{m}}s_m +c_{i{(m+1)}}s_{(m+1)}   \in \Ideal{\pset}  \]}

	\noindent
	By induction hypothesis, we have {$ b_{2{m}}^{\EXP(m)}c_{2{m}}f_{m}+c_{2{m}}s_{m}  \in \Ideal{\pset}$} and hence 
	\[ \p_4 =  q_{2(m+1)}c_{i{(m+1)}}b_{2{m}}^{\EXP(m)}b_{1{m}}^{(\EXP(m)-1)} c_{2{m}} f_{m} + c_{i{(m+1)}}s_{(m+1)} \in \Ideal{\pset}  \]
	
	\noindent
	Using { $ q_{2(m+1)}b_{3m}b_{i(m+1)}c_{i(m+1)}f_{m} + q_{2(m+1)}b_{2m}c_{i(m+1)}f_{m} \in \pset$}, $\EXP(m)$ times, we get 
	\[  \p_5 = q_{2(m+1)} b_{i{(m+1)}}^{\EXP(m)}c_{i{(m+1)}}b_{3{m}}^{\EXP(m)}b_{1{m}}^{(\EXP(m)-1)} c_{2{m}}f_{m} + c_{i{(m+1)}}s_{(m+1)} \in \Ideal{\pset} \]
	\noindent
	As { $ { q_{3(m+1)} c_{3m} f_{m} + q_{2(m+1)}}c_{2m}f_{m} \in \pset$,} we obtain
	{\[ \p_6 =  q_{3(m+1)} b_{i{(m+1)}}^{\EXP(m)} c_{i{(m+1)}} b_{3{m}}^{\EXP(m)}  b_{1{m}}^{(\EXP(m)-1)} c_{3{m}} f_{m}+ c_{i{(m+1)}} s_{(m+1)}  \in \Ideal{\pset}  \]}
	
	\noindent
	By induction hypothesis, we have { $ b_{3{m}}^{\EXP(m)}c_{3{m}} f_{m} + c_{3m} s_m  \in \Ideal{\pset}$} and hence 
	\[ \p_7 =  q_{3(m+1)} b_{i{(m+1)}}^{\EXP(m)} c_{i{(m+1)}}  b_{1{m}}^{(\EXP(m)-1)} c_{3{m}}   s_{m} + c_{i{(m+1)}} s_{(m+1)} \in \Ideal{\pset} \]
	
	\noindent
	Since {$q_{3(m+1)}b_{1m}c_{3m}s_{m} + q_{2(m+1)}b_{4m}c_{2m}s_{m} \in \pset$},  we have
	\[  \p_8 =  q_{2(m+1)} b_{i{(m+1)}}^{\EXP(m)} c_{i{(m+1)}}  b_{4{m}}b_{1{m}}^{(\EXP(m)-2)}c_{2{m}}    s_{m} + c_{i{(m+1)}} s_{(m+1)} \in \Ideal{\pset} \]
	\noindent
	Notice that $\p_3=c_{i(m+1)}s_{(m+1)}+b_{1m}\alpha$ and $\p_8=c_{i(m+1)}s_{(m+1)}+b_{4m}b_{i(m+1)}^{\EXP(m)}\alpha$, where $\alpha=q_{2(m+1)} c_{i{(m+1)}} b_{1{m}}^{(\EXP(m)-2)}c_{2{m}}    s_{m}$. That is, the effect of the steps involved, replaces $b_{1m}$ by $b_{4{m}}b_{i{(m+1)}}^{\EXP(m)}$. Iterating this procedure $\EXP(m-1)-2$ times yields

	$$ \p_9=q_{2(m+1)} b_{i(m+1)}^{\EXP(m)(\EXP(m)-1)} c_{i(m+1)}  b_{4m}^{(\EXP(m)-1)} c_{2m} s_m + c_{i{(m+1)}}s_{(m+1)}  \in \Ideal{\pset} $$
	Again, repeating the steps involved in showing $\p_3\in 
	\Ideal{\pset}$ to $\p_7\in \Ideal{\pset}$, we get	
	$$ \p_{10}= q_{3(m+1)} b_{i{(m+1)}}^{\EXP(m) \times \EXP(m)} c_{i{(m+1)}}  b_{4{m}}^{(\EXP(m)-1)} c_{3{m}} s_{m} + c_{i{(m+1)}}s_{(m+1)} \in \Ideal{\pset} $$
	
	\noindent
	Note that $b_{i{(m+1)}}^{\EXP(m) \times \EXP(m)} = b_{i{(m+1)}}^{\EXP(m+1)}$. Now using $ q_{4(m+1)}b_{4m}c_{4m}f_{m} + q_{3(m+1)}c_{3m}s_{m}  \in \pset$ yields  
	\[ \p_{11} = q_{4(m+1)} b_{i{(m+1)}}^{\EXP(m+1)} c_{i{(m+1)}}  b_{4{m}}^{\EXP(m)} c_{4{m}} f_{m} + c_{i{(m+1)}}s_{(m+1)}  \in \Ideal{\pset} \]
	
	\noindent
	By induction hypothesis, we have that $ b_{4{m}}^{\EXP(m)} c_{4{m}}f_{m} + c_{4m}s_m \in \Ideal{\pset}$,
	\[ \p_{12} = q_{4(m+1)} b_{i{(m+1)}}^{\EXP(m+1)} c_{i{(m+1)}}  c_{4{m}} s_{m} + c_{i{(m+1)}} s_{(m+1)} \in \Ideal{\pset} \]
	
	\noindent
	Now using { $q_{4(m+1)}c_{4m}s_{m} + f_{(m+1)} \in \pset $}, we get
	\[ b_{i{(m+1)}}^{\EXP(m+1)} c_{i{(m+1)}} f_{(m+1)} + c_{i{(m+1)}}s_{(m+1)}\in \Ideal{\pset}  \]
\end{proof}

\subsection{Proof of Lemma \ref{lem:mayrmeyer}}

Lemma \ref{lem:mayrmeyer} states that for any $ j \in [1,4]$ and $\alpha \in \monoRes_\pset(c_{jn} s_n)$,  if $\divides{f_n} {\alpha}$ then $\alpha = b_{jn}^{\EXP(n)}c_{jn}f_n $. 

For any $i \leq n$, let $\Pset_i = \bigcup_{l =0}^i \pset_l$ and $\alpha \in   \monoRes_{\Pset_i}(c_{ji}s_i)$. By Lemma \ref{lem:seq}, we have the monomial sequence $\monomial_0, \ldots, \monomial_r$ such that $\monomial_0 = c_{ji}s_i$, $\monomial_r = \alpha$ and $\monomial_k + \monomial_{k+1} = \beta_k \poly_k$, where $\beta_k\in \varset^\oplus$, and $\poly_k\in \Pset_i$.
We define $\height(\monomial) = \texttt{min}(\{ i \in[1,n] \mid \exists j \in [1,4], \divides{c_{ji}}{\monomial}\})$. For example, $\height(q_{1{(m+1)}}c_{2{(m+1)}}c_{1{m}} s_{m} )$ is $m$. 
Now, we have the following observations.
\begin{enumerate}
		\item For all $\alpha\in \monoRes_{\Pset_{i}}(c_{ji}s_i),\ \height(\alpha)$ is well defined. This follows from a simple observation on $\pset$ that $c_{ji}$ divides $\alpha$ for some $j\in[1,4]$.
		\item For any $\mathfrak{j} \in [0,r-1]$, $|\height(\monomial_{\mathfrak{j}+1}) - \height(\monomial_\mathfrak{j})| \leq 1$.  That is, a single step reduction increases or decreases the height by at most $1$. \label{obs:two}
		\item If for any $\mathfrak{j} \in [0,r-1]$, $\height(\monomial_{\mathfrak{j}+1}) = \height(\monomial_\mathfrak{j}) $ then $\poly_\mathfrak{j}$ is one of the polynomials of Type (\ref{pm:3}-\ref{pm:6}) or (\ref{pm:8}). That is, only these polynomials preserve the height.
		\item If for any $\mathfrak{j} \in [0,r-1]$, $\height(\monomial_{\mathfrak{j}+1}) \neq \height(\monomial_\mathfrak{j}) $ then $\poly_\mathfrak{j}$ is one of the polynomials of Type (\ref{pm:2}) or (\ref{pm:7}). Only these polynomials change the height. 
\end{enumerate} 
The proof we present is similar to the proof in \cite[lemma 7]{MM1982}.
We first prove the following auxiliary Lemma \ref{SC:div}, which we will use in the proof.

\begin{lemma}\label{SC:div}
	For any $j \in [1,4]$, let $\monomial \in \monoRes_\pset(c_{jn}s_n)$ and let $k = \height(\monomial)$, then the following are true.
	
	\begin{enumerate}
		\item For all $k \leq m \leq n$, there is a unique $j \in [1,4]$ such that $\divides{c_{j  m}} {\alpha}$.
		\item For all $k \leq m < n$, there is a unique $j \in [1,4]$ such that $\divides{q_{j m}} {\alpha}$.
		\item Either  $\divides{s_k}{\monomial}$ or $\divides{f_k}{\monomial}$, moreover, $\notdivides{s_k f_k}{\monomial}$ and for all $m \neq k$, $\notdivides{s_m}{\monomial}$ and $\notdivides{f_m}{\monomial}$. That is, exactly one of $s_k$ or $f_k$ divides $\alpha$.
	\end{enumerate} 
	
\end{lemma}

\begin{proof}
	Since $\monomial \in \monoRes_\pset(c_{jn}s_n)$, 
	by Lemma \ref{lem:seq}, we can find a sequence of monomials $\monomial_0, \ldots, \monomial_r $ such that $\monomial_0 = c_{jn}s_n$, $\monomial_r = \monomial$ and $ \forall i \in [0, r -1],\ \monomial_i + \monomial_{i+1} =  \beta_i \poly_i$,  where $\beta_i \in \varset^\oplus$ and $\poly_i \in \pset$.
	We will inductively show that the three properties of the lemma are true for each of the $\monomial_i$ and hence conclude that they are true for $\monomial$ as well.

\noindent
{\bf Base case:} When $i=0$, notice that $1$ and $3$ are immediate, and $2$ is vacuously true since $\height(\monomial_0) = n$.
	
	\noindent
{\bf Inductive case:}
{Assuming the result is true for $i$, we will show it for $i+1$}. Let $\monomial_i + \monomial_{i+1} = \beta_i \poly_i$ and let $\height(\monomial_i) = k_i$.  Inductively, we assume that the properties are true for $\monomial_i$.

Suppose $\poly_i$ is a polynomial of Type \eqref{pm:1} in $\pset$, then $\poly_i$ is of the form $b_{j0}^2c_{j0}f_0+c_{j0}s_0$ for some $j\in[1,4]$. Then either $\alpha_i=\beta_ib_{j0}^2c_{j0}f_0$ and $\alpha_{i+1}=\beta_ic_{j0}s_0$ or vice versa.  Then Property $3$ follows immediately (since $\height(\monomial_i) = 0$, and by the induction hypothesis, the property holds for $\alpha_i$).
 Property $1$ and $2$ also follow from the induction hypothesis and by observing $\poly_i$.
	
	Suppose $\poly_i$ is a polynomial of Type \eqref{pm:2} in $\pset$, then we know that $\divides{s_{k_i}} {\monomial_i}$. Then either $\poly_i= {q_{1{k_i}} c_{1{({k_i}-1)}}} s_{{k_i}-1} + s_{k_i} $ or $ \poly_i= q_{1{(k_i+1)}} c_{1{k_i}} s_{k_i} + s_{k_i+1}$. 
	
	If $\poly_i=q_{1{k_i}} c_{1{({k_i}-1)}} s_{{k_i}-1}+s_{k_i}$, then $\monomial_i + \monomial_{i+1} = \beta_i ({q_{1{k_i}} c_{1{({k_i}-1})}} s_{{k_i}-1} + s_{k_i}) $. From this, we can see that $\monomial_i=\beta_i s_{k_i}$. Then we have $\monomial_{i+1} = \frac{\monomial_i}{s_{k_i}} \times {q_{1{k_i}} c_{1{({k_i}-1)}}s_{{k_i}-1} }$. Since we are using Type (2) polynomial, $k_i$ cannot be $0$. Since $\divides{c_{1(k_i-1)}}{\monomial_{i+1}}$ (and $\height(\monomial_i) = k_i$) we also have $\height(\monomial_{i+1}) = {k_i}-1$. 
	Inductively we have for any $k_i \leq m \leq n$, there is a unique $l \in [1,4]$ such that $\divides{c_{l m}} {\monomial_i}$, hence $\divides{c_{l m}} {\monomial_{i+1}}$ as well. In addition, we have $\divides{c_{1(k_i-1)}}{\monomial_{i+1}}$, with this Property $1$ is satisfied. Reasoning for Property $2$ is  similar.  Property $3$ is immediate from the fact that $\monomial_{i+1} = \frac{\monomial_i}{s_{k_i}} \times {q_{1(k_i)} c_{1{(k_i-1)}} s_{k_i-1} }$.
	
	If $\poly_i= {q_{1{(k_i+1)}} c_{1{k_i}} s_{k_i}}  + s_{k_i+1}$, then $\monomial_i + \monomial_{i+1} = \beta_i ({q_{1{(k_i+1)}} c_{1{k_i}} s_{k_i} } + s_{k_i+1} ) $. Since $\height(\alpha_i)=k_i$, we have $\monomial_{i+1}= \frac{\monomial_i}{ {q_{1{(k_i+1)}}c_{1{k_i}}s_{k_i} }} \times s_{k_i+1}$. To see why Property $3$ holds, notice that  $\height(\monomial_{i+1}) = k_i+1$ and that $\divides{s_{k_i+1}}{\monomial_{i+1}}$ (further by induction hypothesis, only $ \divides{s_{k_i}}{\monomial_i}$). The Properties $1,2$ are immediate from the induction hypothesis.
	
	Similarly, suppose  $\poly_i$ is a polynomial of Type ($7$) in $\pset$, then $ \monomial_{i+1}=\frac{\monomial_i}{f_{k_i}} \times {q_{4{k_i}}c_{4{(k_i-1)}}s_{k_i-1}}$ or $\monomial_{i+1} = \frac{\monomial_i}{q_{4{(k_i+1)}} c_{4{k_i}}s_{k_i}} \times f_{k_i+1}$. We can prove this case using the similar arguments as above.  
	
	Suppose $\poly_i$ is a  polynomial of Type ($3$), then $\poly_i = q_{2(k_i+1)}c_{2k_i}s_{k_i} + q_{1(k_i+1)}b_{1k_i}c_{1k_i}f_{k_i}$.  Clearly $\alpha_{i+1} = \frac{\alpha_i}{q_{2(k_i+1)}c_{2k_i}s_{k_i}} q_{1(k_i+1)}b_{1k_i}c_{1k_i}f_{k_i}$ or $\alpha_{i+1} = \frac{\alpha_i}{q_{1(k_i+1)}b_{1k_i}c_{1k_i}f_{k_i}} {q_{2(k_i+1)}c_{2k_i}s_{k_i}}$.  Now by combining this with the induction hypothesis for $\alpha_i$,  it is easy to see that the Properties $1$-$3$ hold for $\alpha_{i+1}$.  The argument is similar when $\poly_i$ is a  polynomial of Type ($4-6$) or ($8$). 
\end{proof}

\noindent
\begin{figure}[h]
	\centering
\includegraphics[trim = 0.5cm 0.5cm 0.5cm 0,scale=0.26]{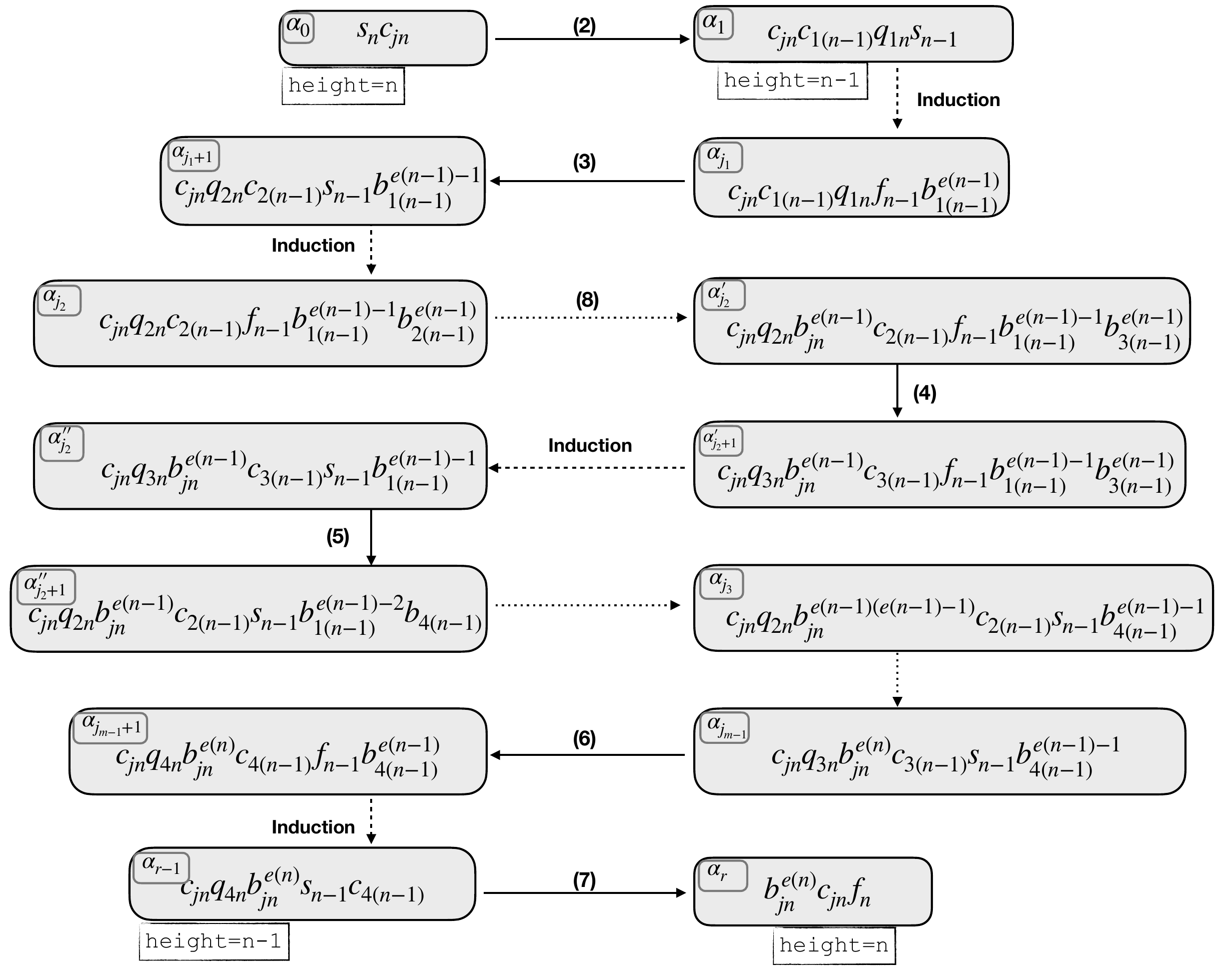}
\caption{Summary of the proof of Proposition \ref{prop:mayr1} and  Lemma \ref{lem:mayrmeyer}}\label{fig:meyermeyer}
\end{figure}

\begin{proof} [Proof of Lemma \ref{lem:mayrmeyer}]
	As in the case of Proposition \ref{prop:mayr1}, here also the Figure \ref{fig:meyermeyer} will provide as a reference.  
	Let $\Pset_i = \bigcup_{l =0}^i \pset_l$, we will show by inducting on $i \in [0,n]$ that for any   $j \in [1,4]$,  and $\beta \in   \monoRes_{\Pset_i}(c_{ji} s_i)$,  if $\divides{f_i} {\beta}$ then $\beta = b_{ji}^{\EXP(i)} c_{ji} f_i $. In fact we will strengthen our induction hypothesis further to include the condition that if $\divides{s_i} {\beta}$, then $\beta =c_{ji} s_i$.  
	
	\paragraph{ Case $i=0$} Recall that $	\pset_0 = \{ b_{j0}^2 c_{j0} f_0+ c_{j0} s_0 \mid j \in [1,4] \}$, notice that for any $j \in [1,4]$, exactly only one polynomial is applicable.	From this the proof is immediate.
	
\paragraph{ Case $i>0$} Let $\beta \in   \monoRes_{\Pset_i}(c_{ji}s_i)$ and let $\divides{f_i}{\beta}$, then by Lemma \ref{lem:seq}, we have the monomial sequence $\monomial_0, \ldots, \monomial_r$ such that $\monomial_0 = c_{ji}s_i$, $\monomial_r = \beta$ and $\monomial_{\mathfrak{j}} + \monomial_{\mathfrak{j}+1} = \beta_\mathfrak{j} \poly_\mathfrak{j}$, where $\beta_\mathfrak{j}\in \varset^\oplus$, and $\poly_\mathfrak{j}\in \Pset_{i}$. We will assume that the sequence is repetition free (i.e. for any $0 \leq \mathfrak{i} < \mathfrak{j}  \leq r$,  $\alpha_\mathfrak{j} \neq \alpha_{\mathfrak{i}}$), this is crucial to our construct.

	Since $\beta \in   \monoRes_{\Pset_i}(c_{ji}s_i) $ and  $\divides{f_i}{\beta}$, by Lemma \ref{SC:div}, we have $\height(\monomial_r)=i$.
	The only possibility for $\poly_0$ is that, it is a polynomial of Type (\ref{pm:2}) in $\pset_i$. From this we will obtain that $\monomial_1 = 
	c_{ji}c_{1{(i-1)}}q_{1{i}}s_{(i-1)} $ and $\height(\monomial_1)<\height(\monomial_0)$.  Let $\mathfrak{j}_1>1 $ be the least index such that $\height(\monomial_{\mathfrak{j}_1  }) = \height(\monomial_1)$.  Note that $\mathfrak{j}_1 < r$, since $\height(\monomial_r) > \height(\monomial_1)$.
	
	{Further for all $\mathfrak{j} \in [2,\mathfrak{j}_1-1]$, $\height(\monomial_\mathfrak{j}) < \height(\monomial_1)$, hence we have that each $\poly_\mathfrak{j} \notin \cup_{l=i-1}^n \pset_l$ (follows from the definition of polynomials). Observe that the variables from $V_i$ do not divide the terms of $\poly_\mathfrak{j}$ (i.e. $\poly_\mathfrak{j}$ does not contain variables from $V_i$)  and $q_{1i} c_{ji}$ divides $\alpha_1$, hence $\divides {q_{1i}c_{ji}}{\alpha_{\mathfrak{j}
+1}}$. From this, we have $\monomial_{\mathfrak{j}_1} =  q_{1{i}}c_{ji}\gamma_{\mathfrak{j}_1}$ for some $\gamma_{\mathfrak{j}_1}$ such that   $\gamma_{\mathfrak{j}_1} \in \monoRes_{\Pset_{i-1}}(c_{1(i-1)}s_{(i-1)})$ and $\height(\gamma_{\mathfrak{j}1})=i-1$. If $\divides{s_{i-1}}{\gamma_{\mathfrak{j}_1}}$ then by induction hypothesis, $\gamma_{\mathfrak{j}1}=c_{1(i-1)}s_{i-1}$, this would contradict repetition free assumption. Then, by Lemma \ref{SC:div}-3, we have $\divides{f_{i-1}}{\gamma_{\mathfrak{j}_1}}$. Inductively we get $\gamma_{\mathfrak{j}_1}=b_{1(i-1)}^{\EXP(i-1)}c_{1(i-1)}f_{(i-1)}$ and hence $\monomial_{\mathfrak{j}_1} =  q_{1{i}}c_{ji} b_{1(i-1)}^{\EXP(i-1)}c_{1(i-1)}f_{(i-1)}$. Now, the only possibility for $\poly_{\mathfrak{j}_1}$ is to be of Type  (\ref{pm:3}) and (\ref{pm:7}). However,  the only repetition free possibility is Type (\ref{pm:3}) from $\pset_{i}$ and this immediately implies that $\monomial_{\mathfrak{j}_1+1} =  q_{2{i}} b_{1(i-1)}^{\EXP(i-1)-1}c_{2(i-1)}c_{ji}s_{(i-1)}$.}
	
	The only polynomials that we can apply to $\monomial_{\mathfrak{j}_1+1}$ are Type (\ref{pm:3}) from $\pset_i$ and Type (\ref{pm:2}) from $\pset_{i-1}$. Again the repetition free possibility is only Type (\ref{pm:2}) from $\pset_{i-1}$.
Then $\monomial_{\mathfrak{j}_1+2} = q_{2{i}} b_{1(i-1)}^{\EXP(i-1)-1}c_{2(i-1)}c_{ji}s_{(i-2)}q_{1(i-1)}c_{1(i-2)}$ and $\height(\monomial_{\mathfrak{j}_1+2}) < \height(\monomial_{\mathfrak{j}_1+1})$. Further, there is an $\monomial_{\mathfrak{j}_2}$ (where $\mathfrak{j}_2$ is the least such index) with $\height(\monomial_{\mathfrak{j}_2}) = \height(\monomial_{\mathfrak{j}_1})$ and $\mathfrak{j}_2 < r$. 
	Since for all $\mathfrak{j} \in [\mathfrak{j}_1+2, \mathfrak{j}_2-1]$, $\poly_\mathfrak{j} \notin  \cup_{l=i-1}^n \pset_l$, we can observe that $\poly_{\mathfrak{j}}$ does not contain the variables $b_{1(i-1)},c_{ji}, c_{2(i-1)}$ and $q_{2i} $. Hence, we have $\divides{b_{1(i-1)}^{\EXP(i-1)-1}q_{2i}c_{2(i-1)}c_{ji}}{\alpha_{\mathfrak{j}+1}}$. In particular, 
	 $\divides{b_{1(i-1)}^{\EXP(i-1)-1}q_{2i}c_{2(i-1)}c_{ji}}{\alpha_{\mathfrak{j}_2}}$.  Now we can rewrite $\alpha_{\mathfrak{j}_2} = q_{2{i}} b_{1(i-1)}^{\EXP(i-1)-1}c_{ji} \gamma_{\mathfrak{j}_2}$, where $\gamma_{\mathfrak{j}_2} \in \monoRes_{\Pset_{i-1}}(c_{2(i-1)}s_{(i-1)}) $. Since     $\height(\alpha_{\mathfrak{j}_2}) = i-1  $, either $\divides{s_{(i-1)}}{\alpha_{\mathfrak{j}_2}} $ or $\divides{f_{(i-1)}}{\alpha_{\mathfrak{j}_2}} $. By the induction hypothesis and the repetition free assumption, we can conclude that $\gamma_{\mathfrak{j}_2} = c_{2(i-1)}f_{(i-1)}b_{2(i-1)}^{\EXP(i-1)}$. Hence we have $\monomial_{\mathfrak{j}_2} =  q_{2{i}} b_{2(i-1)}^{\EXP(i-1)}b_{1(i-1)}^{\EXP(i-1)-1} c_{2(i-1)}c_{ji}f_{(i-1)}$. 
	Now notice that the only possibilities for $\poly_{\mathfrak{j}_2}$ is a polynomial of Type (\ref{pm:4}) and Type (\ref{pm:8}) in $\pset_i$. Suppose we apply  Type (\ref{pm:4}), then the repetition free possibility is only Type (\ref{pm:7}). However notice that $\notdivides{b_{3(i-1)}^{\EXP(i-1)}}{\monomial_{\mathfrak{j2}}}$ and by the induction hypothesis this is required to come back to height $i-1$ (in such a way that the repetition free assumption is not violated). Hence we apply Type (\ref{pm:8}) and this in effect replaces   $b_{2(i-1)}$ by $b_{3(i-1)}$.
	On repeating the argument $\EXP(i-1)$ times, we will obtain $\monomial_{\mathfrak{j}'_2} = q_{2{i}}b_{3(i-1)}^{\EXP(i-1)}b_{ji}^{\EXP(i-1)} b_{1(i-1)}^{\EXP(i-1)-1}c_{2(i-1)}c_{ji} f_{(i-1)}$.
	
	At this point the unique repetition free possibility  for $\poly_{\mathfrak{j}'_2}$ is a polynomial of Type (\ref{pm:4}) of $\pset_i$.
	Hence we have $\monomial_{\mathfrak{j}'_2+1} = q_{3{i}}b_{3(i-1)}^{\EXP(i-1)}b_{ji}^{\EXP(i-1)}b_{1(i-1)}^{\EXP(i-1)-1}c_{3(i-1)}c_{ji}f_{(i-1)}$. 
	   The only possibility for $\poly_{\mathfrak{j}'_2+1} $ is a polynomial of Type (\ref{pm:7}) in $\pset_{i-1}.$ Then $\height(\monomial_{\mathfrak{j}'_2+2})<\height( \monomial_{\mathfrak{j}'_2+1})$. Let $\mathfrak{j}''_2>  \mathfrak{j}'_2+1$ be the least index such that $\height( \monomial_{\mathfrak{j}''_2})=\height( \monomial_{\mathfrak{j}'_2+1})$.
	Then, using an argument similar to the one above, we can inductively obtain  $\monomial_{\mathfrak{j}''_2} =  q_{3{i}}b_{ji}^{\EXP(i-1)}b_{1(i-1)}^{\EXP(i-1)-1}c_{3(i-1)}c_{ji} s_{(i-1)}$.
	We now note that the only possibilities for $\poly_{\mathfrak{j}''_2}$ are polynomials of Type (\ref{pm:5}) and Type (\ref{pm:6}). Suppose we apply Type (\ref{pm:6}), then the repetition free possibility is to apply Type (\ref{pm:7}) of $\pset_{i-1}$. In order to use induction hypothesis, we need the resulting monomial to be divisible by ${b_{4(i-1)}^{\EXP(i-1)} }$ (otherwise we may not be able to return to height $i$), this is not the case.  
	So, we apply Type (\ref{pm:5}), and as a consequence, we have $ \monomial_{\mathfrak{j}''_2+1} = q_{2{i}}b_{4{(i-1)}}b_{ji}^{\EXP(i-1)} b_{1(i-1)}^{\EXP(i-1)-2}c_{2(i-1)}c_{ji}s_{(i-1)}$. 
	Notice that we are in a case similar to $\monomial_{\mathfrak{j}_1+1}$, and we can now repeat the arguments applied so far to get a repetition free path as follows.

	{ The construction from $\monomial_{\mathfrak{j}_1+1}$ to $\monomial_{\mathfrak{j}''_2+1}$ involves replacing $b_{1(i-1)}$ by   $b_{4(i-1)}b_{ji}^{\EXP(i-1)}$. Iterating this procedure ${\EXP({i-1})-2}$ times yields $\monomial_{\mathfrak{j}_3}=q_{2i}c_{ji}b_{4(i-1)}^{\EXP{(i-1)}-1}c_{2(i-1)}s_{i-1}b_{ji}^{\EXP(i-1)(\EXP(i-1)-1)}$. Again, continuing with the same arguments from $\monomial_{\mathfrak{j}_1+1}$ to $\monomial_{\mathfrak{j}_2''}$}, we obtain $\monomial_{\mathfrak{j}_{m-1}} = q_{3{i}}b_{4(i-1)}^{\EXP(i-1)-1}b_{ji}^{\EXP(i)}c_{3(i-1)}c_{ji} s_{(i-1)}$. 
The only possibility for $\poly_{\mathfrak{j}_{m-1}} $ is of Type (\ref{pm:6}) in $\pset_i$. With this, we also obtain $\monomial_{\mathfrak{j}_{m-1}+1} = q_{4{i}}b_{4(i-1)}^{\EXP(i-1)}b_{ji}^{\EXP(i)} c_{4(i-1)}c_{ji}f_{(i-1)}$. Finally applying the induction hypothesis and polynomial of Type (\ref{pm:7}), we obtain
	$\monomial_r = b_{ji}^{\EXP(i)}c_{ji}f_i $ as required.
\end{proof}

\section{Conclusion}

In this work, we considered the question of whether there are polynomials for which the corresponding reduced Gr\"obner basis is  double exponential in size of the given polynomials for different monomial orderings. We showed that there is a family of such polynomials for which the reduced Gr\"obner basis under lexicographic, degree lexicographic and the weighted ordering is double exponential in size. In fact, we identified a sufficient criterion on the monomial ordering, that guarantees a double exponential sized reduced Gr\"obner basis. We believe that our work will allow for the identification of sub-class of polynomials for which computing the  reduced Gr\"obner basis is easy. As a future work, it will be interesting to find whether lexicographic ordering is always harder than degree lexicographic ordering. Another interesting question is to identify conditions under which one monomial ordering performs better than another.
\section*{Acknowledgements}

The second author is thankful to UGC for the financial assistance in the form of Senior Research Fellowship (SRF) and also thankful to the Institute of Mathematical Sciences (IMSc), Chennai for their warm support during the preparation of this paper.  All the authors would like to acknowledge the  grant UOH-IOE-RC5-22003 from the IOE-Directorate, University of Hyderabad.



\begin{thebibliography}{013}
\bibitem[AL1994]{AL1994} William W. Adams, Philippe Loustaunau \emph{An introduction to Gr\"obner Bases}, American Mathematical Society(1994), {https://doi.org/10.1090/gsm/003}

\bibitem[Arn2000]{Arn2000} Elizabeth A. Arnold \emph{Computing Gr\"obner bases with Hilbert lucky primes}, Ph.D. Dissertation, University of Maryland, College Park(2000).

\bibitem[Arn2003]{Arn2003} Elizabeth A. Arnold, \emph{Modular algorithms for computing Gr\"obner bases}, Journal of Symbolic Computation(2003), pp. 403-419, {https://doi.org/10.1016/S0747-7171(02)00140-2}.

\bibitem[BS1988]{BS1988} D. Bayer, M. Stillman, \emph{On the complexity of computing syzygies}, Journal of Symbolic Computation(1988), pp. 135-147, {https://doi.org/10.1016/S0747-7171(88)80039-7}.

\bibitem[Buc2006]{Buc2006} B. Buchberger, \emph{Bruno Buchberger’s PhD thesis 1965: An algorithm for finding the basis elements of the residue class ring of a zero dimensional polynomial ideal}, Journal of Symbolic Computation(2006), pp. 475-511, {https://doi.org/10.1016/j.jsc.2005.09.007}.

\bibitem[CLO1997]{CLO1997} D.Cox, J.Little, D.O'shea, \emph{Ideal, Varieties and Algorithms. An introduction to computational algebraic geometry and commutative algebra}, Springer Science and Business media(1997), {https://doi.org/10.1007/978-0-387-35651-8}.



\bibitem[Dub1990]{Dub1990} T.W. Dube, \emph{The structure of polynomial ideals and Gr\"obner bases}, SIAM Journal on Computing{1990}, pp. 753-773, {https://doi.org/10.1137/0219053}.

\bibitem[Fau1999]{Fau1999} Jean Charles Faug\'{e}re, \emph{A new efficient algorithm for computing Gr\"{o}bner bases {$(F_4)$}}, Journal of Pure and Applied Algebra(1999), pp. 61-88, {https://doi.org/10.1016/S0022-4049(99)00005-5}.

\bibitem[Fau2002]{Fau2002} Jean Charles Faug\'{e}re, \emph{A new efficient algorithm for computing Gr\"{o}bner bases without reduction to zero {$(F_5)$}}, Proceedings of the 2002 International Symposium on  Symbolic and Algebraic Computation(2002), pp. 75-83, {https://doi.org/10.1145/780506.780516}.



\bibitem[Huy1986]{Huy1986} Dung T. Huynh, \emph{A superexponential lower bound for Gr\"obner bases and Church-Rosser commutative thue systems}, Information and Control(1986), pp. 196–206, {https://doi.org/10.1016/S0019-9958(86)80035-3}.

\bibitem[Kal2001]{Kal2001} K. Kalorkoti, \emph{Counting and Gr{\"o}bner bases}, Journal of Symbolic Computation(2001), pp.307-313, {https://doi.org/10.1006/jsco.2000.1575}.

\bibitem[Laz1983]{Laz1983} D. Lazard, \emph{Gr\"obner bases, Gaussian elimination and resolution of systems of algebraic equations(1983)}, pp.146-156. {https://doi.org/10.1007/3-540-12868-9\_99}. 

\bibitem[MM1982]{MM1982} Ernst W. Mayr, Albert R. Meyer, \emph{The complexity of the word problems for the commutative semigrous and polynomial ideals}, Advances in Mathematics(1982), pp. 305-329, {https://doi.org/10.1016/0001-8708(82)90048-2}.

\bibitem[MR2013]{MR2013} Ernst W. Mayr and Stephan Ritscher, \emph{Dimension-dependent bounds for Gr\"obner bases of polynomial ideals}, Journal of Symbolic Computation, pp 78-94, {\bf 49} (2013), {https://doi.org/10.1016/j.jsc.2011.12.018}.

\bibitem[Win1988]{Win1988} Franz Winkler, \emph{A $p$-adic approach to the computation of Gr\"obner bases}, Journal of Symbolic Computation, pp. 287-304, {https://doi.org/10.1016/S0747-7171(88)80049-X}

\bibitem[Yap1991]{Yap1991} Chee K. Yap, \emph{A new lower bound construction for commutative Thue systems with applications}, Journal of Symbolic computation(1991), pp. 1-27, {https://doi.org/10.1016/S0747-7171(08)80138-1}.
\end{thebibliography}
\end{document}